\documentclass[11pt]{amsart}
\usepackage{amsfonts}
\usepackage{amsmath, amscd, amssymb, yhmath}
\usepackage{latexsym}
\usepackage{mathrsfs}
\usepackage{xy}
\input xy
\xyoption{all} 

\def \R{\mathbb{R}}
\def \Z{\mathbb{Z}}
\def \C{\mathbb{C}}

\def \c{\mathbf{C}}
\def \m{\Omega^{\infty}\mathbf{hV}}
\def \I{\mathfrak{I}}
\def \J{\mathcal{J}}

\def \BSO{\mathop{\mathrm{BSO}}}
\def \SO{\mathop{\mathrm{SO}}}

\def \id{\mathrm{id}}
\def \Diff{\mathop{\mathrm{Diff}}}
\def \BDiff{\mathop{\mathrm{BDiff}}}

\def \E{\mathbf{E}}
\def \D{\mathbf{T}}

\newtheorem{theorem}{Theorem}[section]

\newtheorem{lemma}[theorem]{Lemma}
\newtheorem{corollary}[theorem]{Corollary}
\newtheorem{proposition}[theorem]{Proposition}

\theoremstyle{remark}
\newtheorem{remark}[theorem]{Remark}

\newtheorem{definition}[theorem]{Definition}

\newtheorem{example}[theorem]{Example}

\begin{document}
%\hspace{-7mm} \texttt{PRELIMINARY VERSION}
%\vskip 5mm

\title{Stable characteristic classes of smooth manifold bundles}
\author{Rustam Sadykov}
\thanks{The author has been supported by a
Postdoctoral Fellowship of Max Planck Institute, Germany}
\date{\today}

%55R40, 58K30

\begin{abstract} 
Characteristic classes of oriented vector bundles can be identified with cohomology classes of the disjoint union $\sqcup \BSO_n$ of classifying spaces of special orthogonal groups $\SO_n$ with $n=0,1,\dots$. A characteristic class is stable if it extends to a cohomology class of a homotopy colimit $\BSO$ of classifying spaces $\BSO_n$. 
%These are known to be Stiefel-Whitney and Pontrjagin classes (Euler classes %are characteristic, but not stable).  

Similarly, characteristic classes of smooth oriented manifold bundles with fibers given by oriented closed smooth manifolds of a fixed dimension $d\ge 0$ can be identified with cohomology classes of the disjoint union $\sqcup \BDiff M$ of classifying spaces of orientation preserving diffeomorphism groups of oriented closed manifolds of dimension $d$. A characteristic class is stable if it extends to a cohomology class of a homotopy colimit of spaces $\BDiff M$. We show that each rational stable characteristic class of oriented manifold bundles of even dimension $d$ is tautological, e.g., if $d=2$, then each rational stable characteristic class is a polynomial in terms of Miller-Morita-Mumford classes.  

%For example, in the case $d=2$ we show that the ring of rational stable 
%characteristic classes is the polynomial ring in terms of %Miller-Morita-Mumford classes. 
 
%We study characteristic classes of smooth fiber bundles with fibers given by %smooth closed manifolds of an arbitrary fixed dimension. A stabilization of a %fiber bundle is defined to be a new fiber bundle obtained by attaching in %consistent way a copy of a handle to each fiber of the bundle. 
%We say that a stabilization is trivial if the bundle of attached handles is %trivial. A characteristic class of fiber bundles is said to be stable %(respectively weakly stable) if it does not change under stabilizations %(respectively trivial stabilizations). For example, higher %Miller-Morita-Mumford classes and higher Franz-Reidemeister torsion invariants
%are weakly stable characteristic classes over an appropriate class of smooth %fiber bundles. In the current note we compute the rings of all stable and %weakly stable characteristic classes of smooth fiber bundles with coefficients %in a field of characteristic zero. In particular, we show that weakly stable %characteristic classes of surface bundles are given by the %Miller-Morita-Mumford classes. 
\end{abstract}

\maketitle

\setcounter{tocdepth}{1}
\tableofcontents

\addcontentsline{toc}{part}{Introduction}
\section{Introduction}

A \emph{characteristic class} of oriented vector bundles is a function $c$ that assigns a cohomology class $c(\xi)\in H^*(X)$ with coefficients in a field $\Bbbk$ to each oriented vector bundle $\xi$ over a topological space $X$. The function $c$ is required to satisfy the condition $c(f^*\xi)=f^*c(\xi)$ for each continuous map $f: Y\to X$ of topological spaces. A characteristic class $c$ is said to be \emph{stable}, if $c(\xi\oplus \varepsilon)=c(\xi)$ for each oriented vector bundle $\xi$ and a trivial line bundle $\varepsilon$ over $X$. 

It is well known that every characteristic class $c$ can be identified with a cohomology class of the disjoint union $\sqcup \BSO_n$ of classifying spaces of special orthogonal groups $\SO_n$ with $n=0,1,\dots$. A characteristic class $c$ is stable if its components $c_n\in H^*(\BSO_n)$ are compatible in the sense that $c$ extends to a cohomology class of the homotopy theoretic union (called \emph{homotopy colimit}) $\BSO$ of classifying spaces $\BSO_n$ taken with respect to inclusions 
\[
{\BSO}_1\longrightarrow {\BSO}_{2}\longrightarrow {\BSO}_3 \longrightarrow \cdots.   
\]
In the present paper we introduce and study stable characteristic classes of  oriented \emph{manifold bundles}. These classes are defined similarly to those of oriented \emph{vector bundles}. 

An oriented {\it manifold bundle} $\xi$ of dimension $d$ over a topological space $X$ is defined to be a locally trivial fiber bundle $E\to X$ with fiber given by an oriented smooth manifold of dimension $d$. It is convenient to include in the definition the assumption that the total space $E$ of the manifold bundle $\xi$ is embedded into $S^{\infty}\times X$ so that the fiber of $\xi$ over a point $x\in X$ is an embedded smooth submanifold of a copy $S^{\infty}\times \{x\}$ of $S^{\infty}$. 
% of an arbitrary but fixed dimension $d$.  
%In the definition of manifold bundles over general topological spaces, %however, there is a subtlety as one needs, for example, to specify what does %it mean for smooth structures on fiber manifolds to vary continuously. 

A {\it characteristic class} $c$ of %dimension $t$ of 
oriented manifold bundles is a rule that associates a cohomology class $c(\xi)\in H^*(X)$ 
to each oriented manifold bundle $\xi$ of an arbitrary but fixed dimension $d$ over a topological space $X$ so that for each continuous map $f: Y\to X$ of base spaces, there is an equality 
\begin{equation}\label{equ:1.2}
c(f^*\xi)=f^*c(\xi)
\end{equation}
of cohomology classes in $H^*(Y)$. In particular $c$ can be identified with a cohomology class of the disjoint union $\sqcup \BDiff M$ of classifying spaces of orientation preserving diffeomorphism groups of oriented closed manifolds of dimension $d$. 

We are interested in so-called {\it stable} characteristic classes of oriented manifold bundles. These are characteristic classes $c\in H^*(\sqcup\BDiff M)$ with components $c_M\in H^*(\BDiff M)$ that are compatible in the sense that  $c$ extends to an appropriately defined homotopy colimit of classifying spaces.

The compatibility condition is motivated by two remarkable facts: a classical result on  characteristic classes of coverings and a recent solution of the standard Mumford conjecture on characteristic classes of oriented surface bundles. 

In the case of (non-oriented) coverings, the diffeomorphism group $\Diff M$ of the fiber coincides with the group $\Sigma_n$ of permutations of $n$ points. The disjoint union of a universal $n$-sheet covering over $\mathrm{B}\Sigma_n$ and the one sheet covering over $\mathrm{B}\Sigma_n$ form an $(n+1)$-sheet covering over $B\Sigma_n$, which is classified by an inclusion $\mathrm{B}\Sigma_n\to \mathrm{B}\Sigma_{n+1}$. A well known Barratt-Priddy-Quillen~\cite{Ba},~\cite{Pr},~\cite{BP} and Segal~\cite{Seg} theorem implies that the homotopy colimit $\mathrm{B}\Sigma_{\infty}$ of these inclusions
\[
   \mathrm{B}\Sigma_0\longrightarrow \mathrm{B}\Sigma_1\longrightarrow
   \mathrm{B}\Sigma_2\longrightarrow \cdots
\]
is homology equivalent to the space $\Omega_0^{\infty}S^{\infty}$ of pointed loops of degree $0$. In particular, cohomology classes of $\mathrm{B}\Sigma_{\infty}$ can be identified with those of $\Omega_0^{\infty}S^{\infty}$. We note that a characteristic class $c\in H^*(\sqcup \mathrm{B}\Sigma_n)$ of (non-oriented) coverings extends to a cohomology class of   
$\mathrm{B}\Sigma_{\infty}$ if and only if for each covering $\xi$ over a topological space $X$,
\begin{equation}\label{equ:1.1}
    c(\xi\sqcup \mathbf{1}_X)=c(\xi),
\end{equation}
where $\mathbf{1}_X$ is the trivial one sheet covering over $X$. 

In the case of oriented surfaces there is no direct way of comparing classifying spaces $\BDiff F_g$ of orientation preserving diffeomorphism groups of different surfaces $F_g$ of genus $g$ as there is no embedding $F_g\subset F_{g+1}$. However, by the Harer-Ivanov stability theorem~\cite{Ha}, \cite{Iv}, in the stable range of $n<g/2-1$ the $n$-th homology groups of $\BDiff F_g$ are isomorphic to those of the classifying space $\BDiff F_{g,1}$, where $F_{g,1}$ is an oriented surface of genus $g$ with one boundary component and $\Diff F_{g,1}$ is the group of orientation preserving diffeomorphisms of $F_{g,1}$ pointwise trivial on the boundary $\partial F_{g,1}$. On the other hand, since every orientation preserving diffeomorphism of $F_{g,1}$ pointwise trivial on $\partial F_{g,1}$ extends trivially to an orientation preserving diffeomorphism 
of $F_{g+1, 1}\supset F_{g, 1}$, there is a sequence of inclusions
\[
   \BDiff F_{0,1}\longrightarrow \BDiff F_{1,1}\longrightarrow \BDiff F_{2,1}\longrightarrow \cdots.
\]
Recently Madsen and Weiss~\cite{MW} (see also an alternative proof by Galatius, Madsen, Tillmann and Weiss in \cite{GMTW}, and Eliashberg, Galatius and Mishachev in \cite{EGM}) proved that its homotopy colimit $\BDiff F_{\infty, 1}$ is homology equivalent to any path component of a relatively simple infinite loop space $\Omega^{\infty}\C P^{\infty}_{-1}$. This implies that the rational cohomology classes of $\BDiff F_{\infty, 1}$ are given by polynomials in terms of so called Miller-Morita-Mumford classes as it has been predicted by the standard Mumford conjecture. 

The reader may easily formulate a compatibility condition for components of a characteristic class $c\in H^*(\sqcup \BDiff F_{g,1})$ similar to condition (\ref{equ:1.1}).

\section{Results}\label{s2}

\subsection{Homotopy colimit of classifying spaces}
%  \addcontentsline{toc}{subsection}{Homotopy colimit of classifying spaces}

In general we may consider oriented compact manifolds $M_{\alpha}$ of a fixed dimension $d\ge 0$ with a (possibly empty) boundary $\partial M_{\alpha}$ parametrized by an orientation reversing diffeomorphism
\[
    w_{\alpha}\colon \partial M_{\alpha} \longrightarrow S^{i-1}\times S^{d-i},  \qquad 1\le i-1\le d-i. 
\]
Again there are inclusions $i^{\alpha}_{\beta}\colon \BDiff M_{\alpha}\to \BDiff M_{\beta}$ whenever $M_{\beta}$ is obtained from $M_{\alpha}$ either by attaching a handle of positive index along $w$, or by attaching a handle of index $0$. The homotopy theoretic union (i.e., homotopy colimit) 
\[
    \mathfrak{H}_1(d)\colon= \mathop{\mathrm{hocolim}}\ \BDiff M_{\alpha}
\] 
with respect to inclusions $i^{\alpha}_{\beta}$ is defined to be the space obtained from $\sqcup \BDiff M_{\alpha}$ by attaching the mapping cylinders of $i^{\alpha}_{\beta}$ along maps identifying copies of $\BDiff M_{\alpha}$ and $\BDiff M_{\beta}$ in the mapping cylinder of $i^{\alpha}_{\beta}$ with the respective copies in the disjoint union $\sqcup \BDiff M_{\alpha}$.  

Let now $c\in H^*(\sqcup \BDiff M)$ be a characteristic class of oriented manifold bundles with fibers given by closed oriented manifolds $M$ of dimension $d$. It is tempting to define $c$ to be stable if it is given by a restriction of a class in $\mathfrak{H}_1(d)$. Such a definition, however, would not be satisfactory. 

\begin{example}\label{ex:1.0}
For $k\ge 1$ let $RW_{2k}$ denote the characteristic class of surface bundles that coincides on sphere bundles with the  Miller-Morita-Mumford class $M_{2k}$ and is trivial on the other surface bundles. This class extends to a cohomology class in $\mathfrak{H}_1(2)$, though intuitively it should not be stable as it is trivial on all ``stabilized" surface bundles (see subsection~\ref{ss:2.2}).   
\end{example}
 
Surprisingly the extraneous classes such as $RW_{2k}$ can be avoided by allowing manifolds $M_{\alpha}$ in the definition of the homotopy colimit to have more than one boundary component. In this paper we consider a homotopy colimit 
\[
    \mathfrak{H}(d)\colon= \mathop{\mathrm{hocolim}}_{\mathfrak{I}}\ \BDiff M_{\alpha}
\]
where the spaces are indexed by a set $\mathfrak{I}$ of pairs $[M_{\alpha}, w_{\alpha}]$ of diffeomorphism types of oriented compact manifolds $M_{\alpha}$ of dimension $d$ with $k\ge 0$ labeled boundary components, and equivalence classes $[w_{\alpha}]$ of orientation reversing parametrizations  
\[
w_{\alpha}\colon \partial M_{\alpha}\longrightarrow \bigsqcup_{i=0}^{\left\lfloor d+1/2\right\rfloor} S^{i-1}\times S^{d-i}
\] 
of the components of the boundary of $M_{\alpha}$ (see section~\ref{s:1.5}).
% a parametrizing diffeomorphism $g$ is said to be equivalent to $f$ if $g$ is %isotopic to the composition of $f$ and a map that switches the factors %$S^{i-1}$ and $S^{d-i}$ of some of the terms $S^{i-1}\times S^{d-i}$. 
If a manifold $M_{\beta}$ is obtained from a manifold $M_{\alpha}$ by attaching handles of the form $D^{i}\times S^{d-i}$ or $S^{i-1}\times D^{d-i+1}$ in an obvious way by means of the parametrization $w_{\alpha}$ of $\partial M_{\alpha}$, then every orientation preserving diffeomorphism of $M_{\alpha}$ trivial on $\partial M_{\alpha}$ extends trivially to an orientation preserving diffeomorphism of $M_{\beta}$. In particular there is an inclusion 
\begin{equation}\label{eq:1.1}
    \BDiff M_{\alpha}\longrightarrow \BDiff M_{\beta}
\end{equation}
of classifying spaces of diffeomorphism groups. The homotopy colimit $\mathfrak{H}(d)$ is formed by means of inclusions (\ref{eq:1.1}), one for each attaching map.  

The definition of $\mathfrak{H}(d)$, however, needs a careful justification as 
\begin{itemize} 
\item there is no canonical choice of a manifold $M_{\alpha}$ in each class $\alpha$, 
\item the spaces $\BDiff M_{\alpha}$ are defined up to weak homotopy equivalence, and 
\item the maps (\ref{eq:1.1}) are defined up to homotopy. 
\end{itemize}

Thus, a priori neither the existence, nor the uniqueness of the space $\mathfrak{H}(d)$ is obvious; we need to show that it is possible to make choices in the construction so that the diagram of all maps (\ref{eq:1.1}) is strictly commutative, and that different choices result in weakly homotopy equivalent spaces.  

\begin{remark}
There is an example of two strictly commutative diagrams of maps $\mathcal{D}_1$ and $\mathcal{D}_2$ and a bijective correspondence between spaces and maps of $\mathcal{D}_1$ and $\mathcal{D}_2$ such that
\begin{itemize} 
\item the spaces in $\mathcal{D}_1$ are the same as the corresponding spaces in $\mathcal{D}_2$, 
\item the maps in $\mathcal{D}_1$ are homotopic to the corresponding maps in $\mathcal{D}_2$, and 
\item the homotopy colimits for $\mathcal{D}_1$ and $\mathcal{D}_2$ are not weakly homotopy equivalent! 
\end{itemize}
\end{remark}

Our first theorem asserts that $\mathfrak{H}(d)$ is well-defined. 

\begin{theorem}\label{th:2.2} The homotopy colimit $\mathfrak{H}(d)$ of classifying spaces of diffeomorphism groups exists and unique.
\end{theorem}

\begin{definition} A characteristic class $c\in H^*(\sqcup \BDiff M)$ of oriented manifold bundles with fibers given by closed oriented manifolds of dimension $d$ is said to be {\it  stable} if it is the restriction of a class in $H^*(\mathfrak{H}(d))$.  
\end{definition}

In the current paper we consider only cohomology groups with coefficients in an arbitrary fixed field $\Bbbk$ of characteristic zero. Our main result, Theorem~\ref{th:1.1}, describes stable characteristic classes in the case of a general fixed even dimension $d\ge 0$ (for the case where $d$ is odd see the paper of Ebert \cite{Eb}).  It turns out that in the case where $d$ is even all stable characteristic classes can be  described in terms of generalized Miller-Morita-Mumford classes, which are cohomology classes of an infinite loop space $\Omega^{\infty}\mathbf{hV}$ defined by Madsen and  Weiss in \cite{MW} (see section~\ref{s:16}). Let $\sqcup\Omega^{\infty}\mathbf{hV}$ denote the disjoint union of copies of a path component of $\Omega^{\infty}\mathbf{hV}$, one term for each cobordism class of oriented closed manifolds of dimension $d$. The geometric meaning of cohomology classes of $\sqcup\Omega^{\infty}\mathbf{hV}$ is well understood (see section~\ref{sec:17}) and it follows that each class in $\sqcup\Omega^{\infty}\mathbf{hV}$ determines a characteristic class which we call {\it tautological}.  

\begin{theorem}\label{th:1.1} Suppose that the dimension $d\ge 0$ is even. Then each stable characteristic class of oriented manifold bundles of dimension $d$ is tautological. 
\end{theorem}

\begin{remark} In the proof of Theorem~\ref{th:1.1} we do not use a Harer-Ivanov stability type theorem; our Theorem~\ref{th:1.1} is true for any even dimension $d\ge 0$, not only $d=2$.  
\end{remark}

In the case $d=0$, the space $\Omega^{\infty}\mathbf{hV}$ coincides with $\Omega^{\infty}S^{\infty}$. By Theorem~\ref{th:1.1}, in this case there are no non-trivial stable characteristic classes with coefficients in a field of characteristic zero. In the case $d=2$ the space $\Omega^{\infty}\mathbf{hV}$ coincides with  the infinite loop spaces $\Omega^{\infty}\C P^{\infty}_{-1}$ of Madsen and Tillmann \cite{MT}. In particular, Theorem~\ref{th:1.1} implies that in the case of surface bundles each stable characteristic class is a polynomial in terms of Miller-Morita-Mumford classes~\cite{Mi}, \cite{Mo1}, \cite{Mo2}. 

In general the space $\mathfrak{H}(d)$ is {\it not} homology equivalent to 
$\sqcup\Omega^{\infty}\mathbf{hV}$. 

According to the Bousfield-Kan cohomology spectral sequence~\cite{DH}, the space $\mathfrak{H}(d)$ may have classes contributed by the right derived functors
\[
   {\lim}^i\, H^*(\BDiff{M})
\]
of the limit functor. In fact, in the case $d=3$ we construct such a phantom class, i.e., a non-trivial class in $\mathfrak{H}(d)$ that is trivial as a stable characteristic class of manifold bundles.

\subsection{Finite stability conditions}\label{ss:2.2}
%The groups $\Diff M$ for different $M$ share common properties and 
We show that the space $\mathfrak{H}(d)$ is in a sense the minimal space containing $\sqcup \BDiff M$ that is suitable for the definition of stable characteristic classes. As it has been mentioned, the attachment of mapping cylinders in the definition of $\mathfrak{H}_1(d)$ creates many extraneous classes. We show that the space $\mathfrak{H}(d)$ can be build from $\mathfrak{H}_1(d)$ in the spirit of the Quillen ``+"-construction eliminating extraneous classes. 
 
In fact for each positive integer $k$ we introduce a stability condition of order $k$ on $c$ so that the order of stability measures the compatibility of components $c_M\in H^*(\BDiff M)$ of $c$. These are defined so that stable characteristic classes are those that satisfy the stability condition of each order $k$.  

Namely let $\I(k)\subset \I$ denote the set of indexes $\alpha$ of diffeomorphism types of oriented compact manifolds $M_{\alpha}$ with at most $k$ parametrized boundary components, each of the form $S^{i-1}\times S^{d-i}$. 

\begin{definition}\label{d:1.1} A characteristic class $c$ is said to be {\it stable of order $k$} if it extends to the homotopy colimit
\[
    \mathfrak{H}_k(d)\colon= \mathop{\mathrm{hocolim}}_{\alpha\in \mathfrak{I}(k)}\ \BDiff M_{\alpha}.
\]
\end{definition}

\begin{example}\label{ex:1.1}
In section~\ref{s:2} we will see that the class $RW_{2k}$ defined in Example~\ref{ex:1.0} is a stable characteristic class of order $1$. On the other hand $RW_{2k}$ is trivial on almost all fiber bundles. As we will see, this implies that $RW_{2k}$ does not satisfy some higher stability condition.  
\end{example}

%To define $k$-th stability condition we use a Vassiliev resolution type %construction for $\mathfrak{H}(d)$ which is equivalent to a filtration 
%\[
%    \sqcup\BDiff M \cong \D_0\subset \D_1\subset \D_2 \subset \cdots  \subset %\D\cong \mathfrak{H}(d)
%\]
%of a CW complex $\D$ weakly homotopy equivalent to the homotopy colimit %$\mathfrak{H}(d)$. The $k$-th stability condition is the obstruction to %extending a characteristic class $c$ to a cohomology class in $\D_k$. 

We give a geometric interpretation of the $k$-th stability condition and prove a criterion for a  stable characteristic class of order $k-1$ to be of order $k$.

%It is worth mentioning that stable characteristic classes can be defined not %only for families of surfaces but also for example for families of surfaces %with Lefschetz singularities. These are related to Deligne-Mumford %compactifications of moduli spaces of Riemann surfaces~\cite{EG}.    

%Our existence and uniqueness Theorem~\ref{th:2.2} is closely related to the subject of the paper \cite{MW} by Madsen %and Weiss. The authors introduce classifying spaces of certain genuine objects and classifying spaces of their %respective homotopy theoretic counterparts. Homotopy theoretic counterparts are simpler and studied by means of %homotopy theory. The study of classifying spaces of genuine objects is more technically demanding and rests on %homotopy colimit decompositions. Our existence and uniqueness Theorem~\ref{th:2.2} is closely related to Madsen-Weiss %study of classifying spaces of genuine objects (see outline of the proof of the Mumford Conjecture \cite[Section %1.3]{MW}, the details in \cite[Chapter 5]{MW}, and a discussion of the mentioned relation in section \ref{MW} below). 

\subsection*{Acknowledgment} I am thankful to Oscar Randal-Williams for comments, in particular for pointing out to the mentioned characteristic classes $RW_{2k}$ (see Example~\ref{ex:1.1}). I am thankful to David Ayala for a discussion of homotopy colimits (in particular, for Example~\ref{ex:3.2}) and for sharing his PhD thesis with me. Remark~\ref{r:4.5} was pointed out to me by Kiyoshi Igusa, which I gratefully acknowledge.

\section{Motivations/Outline}\label{s:mo}

In this section we outline the content of the paper and shortly explain a few essential points behind our choices made in the paper. 

\subsection{Definition of $\mathfrak{H}(d)$}
\subsubsection{Naive construction} Let $M_{\alpha}$ be a smooth compact oriented manifold with one 
(possibly non-connected) component boundary $\partial M$, which is parametrized by an orientation reversing diffeomorphism 
$w\colon \partial M_{\alpha} \to S^{i-1}\times S^{d-i}$ for some non-negative $i$ with $i-1\le d-i$. Then $M_{\alpha}$ is a submanifold of oriented closed manifolds
\[
   M_0=M_{\alpha}\bigcup_{w} D^i\times S^{d-i}\quad \mathrm{and}\quad M_1=M_{\alpha}\bigcup_{w} S^{i-1}\times D^{d-i+1}
\]
where $w$ identifies $\partial M_{\alpha}$ with the boundary of the corresponding handle. Each orientation preserving diffeomorphism of $M_{\alpha}$ pointwise trivial on $\partial M_{\alpha}$ extends trivially to orientation preserving diffeomorphisms of $M_0$ and $M_1$. In particular there is a diagram of injective homomorphisms
\[
 \Diff M_0 \longleftarrow \Diff M_{\alpha} \longrightarrow \Diff M_1,
\]
and therefore, a diagram of inclusions of classifying spaces 
\begin{equation}\label{eq:3.2}
\BDiff M_0    \stackrel{i_{\alpha}}\longleftarrow \BDiff M_{\alpha} \stackrel{j_{\alpha}}\longrightarrow   \BDiff M_1.    
\end{equation}
The weakest compatibility condition on components $c_M\in H^*(\BDiff M)$ of a characteristic class $c$ requires that $i^*_{\alpha}c_{M_0}=j^*_{\alpha}c_{M_{1}}$ for each diffeomorphism type $\alpha$ of oriented compact manifolds with one boundary component. This compatibility condition is equivalent to the requirement that the components $c_{M_0}$ and $c_{M_1}$ extend to a cohomology class of the colimit of the diagram (\ref{eq:3.2}), which is the space
\[
    \BDiff M_0 \bigcup_{\BDiff M_{\alpha}} \BDiff M_1.
\]
It is tempting to replace the infinitely many conditions, one for each $\alpha$, by a single condition; namely, to  define a stable characteristic class to be the restriction $c \in H^*(\sqcup\BDiff M)$ of a cohomology class of the colimit 
\[
   \mathop{\mathrm{colim}} \BDiff M\colon= \bigsqcup \BDiff M / \sim
\] 
with respect to the diagrams (\ref{eq:3.2}). Here the disjoint union ranges over all diffeomorphism types of oriented closed manifolds $M$ of dimension $d$; and for each equivalence type $[M_{\alpha}, w_{\alpha}]$ of an oriented compact manifold with one parametrized boundary component, the equivalence $\sim$ identifies each point $i_{\alpha}(x)$ with the point $j_{\alpha}(x)$ as above. 

Such a definition, however, would not be satisfactory since  
\begin{itemize}
\item colimits do not behave well under homotopy equivalences, and
\item such a colimit has many extraneous classes.   
\end{itemize}

\subsubsection{Colimits vs Homotopy colimits}
Colimits do not behave well under homotopy equivalences (e.g., see \cite{DS}). For example, if $S^{n-1}\subset D^n$ is the inclusion of the boundary, then the colimits of the upper and lower rows of the diagram
\[
\begin{CD}
  D^n @<\supset<< S^{n-1} @>\subset>> D^n \\
  @VVV @VVV @VVV \\
  pt @<<< S^{n-1} @>>> pt
\end{CD}
\] 
are not homotopy equivalent. To rectify the definition, we replace colimits by homotopy theoretic unions called homotopy colimits. 

\subsubsection{Extraneous classes} The homotopy colimit in the naive definition has many extraneous classes such as $RW_{2k}$. To avoid these classes we consider not the union of diagrams (\ref{eq:3.2}), but its more elaborate version that includes classifying spaces $\BDiff M$ of orientation preserving diffeomorphisms of oriented compact manifolds $M$ with boundary $\partial M$ whose boundary $\partial M$ consists of more than one component, each diffeomorphic to $S^{i-1}\times S^{d-i+1}$. 

\subsection{The Existence Theorem} In sections~\ref{s:1.5}--\ref{s:6} we give a careful definition of $\mathfrak{H}(d)$. One of the main points here is that the homotopy colimit of a diagram is well defined/exists only if the diagram is strictly commutative. Thus, in order to define $\mathfrak{H}(d)$, we need to carefully choose a specific model for each space $\BDiff M_{\alpha}$ and a specific model for each map $\BDiff M_{\alpha}\to \BDiff M_{\beta}$ so that the diagram of maps is strictly commutative. 

\subsection{The Uniqueness Theorem} We warn the reader that, a priori, for a different choice of models for spaces $\BDiff M_{\alpha}$ and maps $\BDiff M_{\alpha}\to \BDiff M_{\beta}$ the corresponding homotopy colimit may not be weakly homotopy equivalent to $\mathfrak{H}(d)$ (see Example~\ref{ex:3.2}). In section~\ref{se:14} we show that the weak homotopy equivalence type of the homotopy colimit $\mathfrak{H}(d)$ does not depend on our choices.

\subsection{Marked fold maps} In order to study stable characteristic classes, we introduce \emph{marked fold maps} (see sections~\ref{se:7}--\ref{se:9}). Their relation to $\mathfrak{H}(d)$ is established later by means of a function $\mathcal{F}$, which we introduce in section~\ref{se:10}. The function $\mathcal{F}$ assigns to each manifold without boundary $N$ the space of marked fold maps to $N$. The function $\mathcal{F}$ resembles a sheaf on the category of smooth manifolds and smooth maps. However a smooth map $N'\to N$ of manifolds does not define a map $\mathcal{F}(N)\to \mathcal{F}(N')$ in general, and, in particular, the definition of a classifying space of a sheaf (see \cite{MW}) is not applicable for the function $\mathcal{F}$. Instead we consider the simplicial set $\mathcal{F}^{\pitchfork}$ of subspaces $\mathcal{F}^{\pitchfork}(\Delta^n_e)\subset \mathcal{F}(\Delta^n_e)$ of those marked fold maps $f$ into extended simplices 
\[
    \Delta^n_e\colon=\{\ (x_0, ..., x_n)\in \R^{n+1}\ |\ x_0+ \cdots +x_n=1\}
\]
that restrict to marked fold maps into subsimplices $\Delta^m_e\subset \Delta^n_e$. We point out that the spaces $\mathcal{F}^{\pitchfork}(\Delta^n_e)$ are essentially different from the spaces $\mathcal{F}(\Delta^n_e)$ (for example, the numbers of path components of these spaces are essentially different in general). Nevertheless, in sections~\ref{se:10}--\ref{se:11} we will show that the realization $|\mathcal{F}|$ of $\mathcal{F}^{\pitchfork}$ yields the space with properties that would be expected for the classifying space of $\mathcal{F}$. In sections~\ref{se:12}--\ref{se:14} we establish the relation of marked fold maps to $\mathfrak{H}(d)$ by showing that $|\mathcal{F}|$ is weekly homotopy equivalent to $\mathfrak{H}(d)$. 

In section~\ref{se:16} we complete the proof of Theorem~\ref{th:1.1}. This part heavily relies on the Eliashberg-Mishachev theorem~\cite{EM} which we recall in section~\ref{s:5}.  In later sections we study the geometric meaning of the finite stability conditions.

\addcontentsline{toc}{part}{Homotopy colimit of classifying spaces}

\section{Manifolds with parametrized boundary}\label{s:1.5}

%To summarize, each choice of an equivalence relation leads to a homotopy %colimit diagram. On the other hand there are several reasonable competing %definitions of equivalence relations that lead to non-equivalent homotopy %colimit diagrams.
%
%\begin{example} Let $M_{\alpha}$ be a disjoint union of two copies of an %orientable closed connected manifold $M$, and $M_{\beta}$ a manifold %diffeomorphic to $M\setminus (S^{1}\times D^{d-1})$. For a given embedding %$i: M_{\beta}\to M_{\alpha}$, let $i'$ denote the embedding given by the one %of the following compositions:
%\begin{itemize}
%\item the composition of $i$ and the permutation of components of %$M_{\alpha}$,
%\item the composition of $i$ and a diffeomorphism of $M_{\beta}$ pointwise %trivial on $\partial M_{\beta}$, 
%\item the composition of $i$ and a diffeomorphism of $M_{\beta}$ %componentwise trivial on $\partial M_{\beta}$, or  
%\item the composition of $i$ and a diffeomorphism of $M_{\beta}$.
%\end{itemize}
%In each case by choosing an equivalence relation we need to decide if $i$ and %$i'$ represent the same arrow in the diagram of the homotopy colimit %$\mathfrak{H}(d)$.
%\end{example}

Let $M$ be an oriented compact (not necessarily path connected) manifold of dimension $d$ with $k\ge 0$ ordered boundary components $\partial_1 M, \partial_2 M, ..., \partial_k M$, each of the form $S^{i-1}\times S^{d-i}$ for some $0\le i\le \frac{d+1}{2}$. 

\begin{remark} Let us recall the convention that for each $d\ge 0$, there is a unique manifold of dimension $d$ whose underlying set is empty. 
\end{remark}

We note that the boundary components of the form $S^{-1}\times S^{d}$ are empty, while the boundary components of the form $S^{0}\times S^{d-1}$ are not path connected. We will assume that the order of boundary components $\partial_1M, \partial_2M, ..., \partial_k M$ agrees with the order given by the indices $i$ of the tori $S^{i-1}\times S^{d-i}$.

\begin{definition}
 A {\it parametrization of the boundary} $\partial M$ is a choice of a map
\[
   w: \partial M\longrightarrow \bigsqcup_{i=0}^{i=\left\lfloor \frac{d+1}{2}\right\rfloor} S^{i-1}\times S^{d-i},
\] 
each of $k$ components of which is an orientation reversing diffeomorphism. We will call a pair $(M, w)$ an \emph{oriented (compact) manifold of dimension $d$ with parametrized boundary}. 
\end{definition}

\begin{remark} The condition that $w$ is orientation reversing guarantees that the orientation of $M$ extends to an orientation of the manifold obtained from $M$ by attaching an oriented handle $D^{i}\times S^{d-i}$ or $S^{i-1}\times D^{d-i+1}$ along any component of $w$. We recall the convention that the boundary $\partial M$ of an oriented manifold $M$ of dimension $d$ is oriented so that an ordered basis $v_1, ..., v_{d-1}$ for the tangent space $T_x\partial M$ at a point $x\in \partial M$ gives an orientation on $\partial M$ at $x$ if the ordered basis $v_0, v_1, ..., v_{d-1}$ for $T_xM$ gives an orientation on $M$ at $x$ where $v_0$ is an arbitrary outward pointing vector at $x$.      
\end{remark}

%\begin{remark} A parametrization $w$ of the boundary $\partial M$ of $M$ is a diffeomorphism onto the image if and %only if $\partial M$ contains at most one path component of each type.  
%\end{remark}

\begin{definition}\label{def:3.1} We say that manifolds with parametrized boundary $(M_1, w_1)$ and $(M_2, w_2)$ are \emph{equivalent} if there is an orientation preserving diffeomorphism 
\[
\alpha: M_1\to M_2
\] 
that preserves the order of boundary components and such that the diagram 
\[
\xymatrix
{
& \partial M_1 \ar[dr]^{w_1} \ar[rr]^{\alpha|\partial M_1}&       & \partial M_2 \ar[dl]_{w_2}     \\     
&            & \bigsqcup_{i} S^{i-1}\times S^{d-i}   & 
}
\]
commutes, where the disjoint union ranges through $i=0, ...,\left\lfloor \frac{d+1}{2}\right\rfloor$.  The equivalence class of a pair $(M, v)$ will be denoted by $[M, v]$.   
\end{definition}

\begin{remark}\label{r:4.3} Let $M_1$ be an oriented closed manifold with no boundary components and let $M_2$ be an oriented closed manifold diffeomorphic to $M_1$  with $k$ ordered boundary components of the form $S^{-1}\times S^d$. The parametrizations $w_1$ of $\partial M_1$ and $w_2$ of $\partial M_2$ are, of course, unique. We note that $(M_1, w_1)$ is not equivalent to $(M_2, w_2)$.  
\end{remark}

%\begin{remark} Let $(M, w)$ be a manifold with parametrized boundary. Suppose $\partial M$ consists of two empty %components $\partial M_1$ and $\partial M_2$. Then the identity map of $M$ is not an equivalence of 
%$(M, \partial_1M, \partial_2M)$ and $(M, \partial_2M, \partial_1M)$. 
%\end{remark}

We note that any two parametrizations $v$ and $w$ of $\partial M$ determine an orientation preserving diffeomorphism $\partial M\to \partial M$ that does not permute path components of $\partial M$. According to Definition~\ref{def:3.1}, $(M, v)$ is equivalent to $(M, w)$ if and only if the corresponding diffeomorphism $\partial M\to \partial M$ extends to a diffeomorphism $M\to M$. 

In particular an \emph{automorphism} of a pair $(M, w)$ is an orientation preserving diffeomorphism of $M$ trivial on $\partial M$. 

\begin{example}\label{ex:3.1} If parametrizations $w_1$ and $w_2$ of the boundary of the same manifold $M$ are isotopic, then $(M, w_1)$ and $(M, w_2)$ are equivalent. 
\end{example}

\begin{remark}\label{r:4.5} The converse of the statement in Example~\ref{ex:3.1} is not true in general. I am thankful to Prof. Igusa for pointing out the following example. Let $\mathcal{P}(M, \partial M)$ be the space of relative pseudo-isotopies of $\mathop{\mathrm{id:}} M\to M$, i.e., diffeomorphisms 
\[
M\times [0,1]\longrightarrow M\times [0,1]
\] 
that keep $M\times \{0\}$ and $\partial M\times [0,1]$ pointwise fixed. Suppose that $\dim M\ge 7$. Then, by a result due to Hatcher~\cite{Hat}, the group 
\[
     \pi_0(\mathcal{P}(M, \partial M))
\]
is trivial if and only if $\pi_1(M)$ is trivial. In particular, for example, there is a diffeomorphism 
\begin{equation}\label{eq:3.1}
\varphi\colon S^1\times S^6\times [0,1] \longrightarrow S^1\times S^6\times [0,1],
\end{equation}
which is trivial on $S^1\times S^6\times \{0\}$ and not isotopic to the identity map. 

We may use the Hatcher Theorem to show that the boundary $\partial M$ of $M\cong S^1\times S^6\times [0,1]$ admits two equivalent non-isotopic parametrizations. Indeed, let $v$ be the parametrization 
\[
    S^1\times S^6\times \{0,1\} \longrightarrow S^1\times S^6
\]
that maps the path component $S^1\times S^6\times \{0\}$ identically on its copy and takes $S^1\times S^6\times \{1\}$ to its copy by flipping the $S^1$ component so that $v$ reverses the orientation of the second path component of $\partial M$.  Let $w$ be the parametrization of $\partial M$ given by $v\circ  \varphi|\partial M$. Then $(M, v)$ is equivalent to $(M, w)$, but $v$ is not isotopic to $w$. 
%%
%
%the parametrization given by the restriction of the map (\ref{eq:3.1}) to the %boundary is equivalent to the identity parametrization, but it is not isotopic %to the identity parametrization. 
\end{remark}

\begin{example} In the case where $d<3$, any two isotopic parametrizations are equivalent.
\end{example}

\section{A model for $\BDiff M$}\label{s:model}

In order to define  the homotopy colimit $\mathfrak{H}(d)$, we will need a specific model for classifying spaces $\BDiff M$, one for each diffeomorphism type of oriented compact manifolds $M$ with parametrized boundary. The model for spaces $\BDiff M$ and maps (\ref{eq:1.1}) should be chosen carefully so that the diagram of all maps (\ref{eq:1.1}) is strictly commutative.   

%To this end let us define a specific model for each space $\BDiff M_{\gamma}$ %and a specific model for each map (\ref{eq:-1}). These specific models for %spaces and maps will be our choices for values of $\mathcal{D}$. 

\begin{remark}\label{re:5.1} The homotopy colimit can be defined only if the diagram of maps (\ref{eq:1.1}) is strictly commutative. Furthermore, the homotopy classes of maps (\ref{eq:1.1}) are canonical, but the maps itself are not canonical. A priori, different choices of maps (\ref{eq:1.1}) in the canonical homotopy classes lead to homotopy colimits that are not weakly homotopy equivalent. We will show, however, that in our particular case, the homotopy colimit does not depend on these choices (see section~\ref{se:14}). 
\end{remark}

The main result of the current section is Theorem~\ref{th:5.1}, whose proof will occupy the rest of the section. 

\begin{theorem}\label{th:5.1}
The models for spaces $\BDiff M$ and maps (\ref{eq:1.1}) can be chosen so that the diagram of all maps (\ref{eq:1.1}) is strictly commutative.
\end{theorem}

\subsection{Almost standard tori and solid tori}
We say that an embedding of a torus $S^{i-1}\times S^{d-i}\to \R^{\infty}$ for some $i=1,..., \left\lfloor \frac{d+1}{2}\right\rfloor$ is {\it almost standard} if it is obtained by the composition of embeddings 
\begin{eqnarray}
    S^{i-1}\times S^{d-i}&\longrightarrow& \R^{i}\times \R^{d-i+1}\times \{pt\}\times \{0\}\times \{0\} \nonumber \\
 &\stackrel{\subset}\longrightarrow& \R^{i}\times \R^{d-i+1}\times \R\times \R\times \R^{\infty-d-3} \nonumber \\
 &\stackrel{\cong}\longrightarrow& \R^{\infty} \nonumber,
\end{eqnarray}
where the first map is the product of two standard embeddings of unit spheres into a copy of $\R^i\times \R^{d-i+1}$ and $pt$ is an arbitrary point in $\R$. Each almost standard torus $S^{i-1}\times S^{d-i}$ bounds two embedded {\it almost standard solid tori}. Namely, 
\[
    D^{i}\times S^{d-i} \subset \R^i\times \R^{d-i+1}\times \{pt\}\times \R\times \{0\} 
\]
given by the set of $5$-tuples $(x_1, x_2, x_3, x_4, x_5)$ with 
\[
    \left\|x_1\right\|\le 1, \quad  \left\|x_2\right\|=1, \quad x_3=pt, \quad x_4\ge 0, \quad \left\|x_1\right\|^2+x_4^2=1, \quad x_5=0; 
\]
and 
\[
    S^{i-1}\times D^{d-i+1} \subset \R^i\times \R^{d-i+1}\times \{pt\}\times \R\times \{0\} 
\]
given by the set of $5$-tuples $(x_1, x_2, x_3, x_4, x_5)$ with 
\[
    \left\|x_1\right\|= 1, \quad  \left\|x_2\right\|\le 1, \quad x_3=pt, \quad x_4\ge 0, \quad \left\|x_2\right\|^2+x_4^2=1, \quad x_5=0. 
\]

\subsection{A model for $\BDiff M_{\alpha}$}\label{ss:5.2} For all $\alpha\in\I$, the space $\BDiff M_{\alpha}$ can be regarded to be the space of equivalence classes of embeddings 
 \[
 f: M_{\alpha}\hookrightarrow \R^{\infty}
 \] 
such that 
\begin{itemize}
\item[(i)] for each non-empty boundary component $\partial_i M_{\alpha}$ with a parametrization $w_i$, the composition $f\circ w_i^{-1}$ is an almost standard embedding of a torus, 
\item[(ii)] the order of the non-empty boundary components of $\partial M$ agrees with the order given by the distances of components of $f(M_{\alpha})$ to $0\in \R^{\infty}$, 
\item[(iii)] the image of $f$ does not intersect the interior of the almost standard solid tori bounded by components of $\partial M$, and
\item[(iv)] if $M_{\beta}$ is obtained from $M_{\alpha}$ by attaching handles along some of the non-empty components of $\partial M_{\alpha}$, then $f$ canonically extends to a smooth embedding $f': M_{\beta}\to \R^{\infty}$ such that $f'|M_{\alpha}=f$ and $f'$ restricted to the glued solid tori coincides with almost standard embeddings of solid tori.   
\end{itemize}

\begin{remark}\ 
\begin{itemize}
\item The condition (iv) implies (iii). 
\item The condition (iv) is included in order to avoid non smooth maps $f'$, e.g., maps onto submanifolds in $\R^{\infty}$ with corners. 
\item The space of embeddings satisfying (i)--(iv) endowed with the $C^{\infty}$ topology is weakly homotopy equivalent to a point. 
\end{itemize}
\end{remark}  
  
An embedding $f$ satisfying (i)--(iv) is said to be {\it equivalent} to an embedding $g$ satisfying (i)--(iv) if there is a diffeomorphism $h$ of $M_{\alpha}$ relative to the boundary such that $f=g\circ h$. In other words, an equivalence class of embeddings satisfying (i)--(iv) is given by a submanifold $W\subset \R^{\infty}$ for which there is a diffeomorphism $f: M_{\alpha}\to W$ as above. The equivalence class of an embedding $f$ is denoted by $[f]$. The space of all equivalence classes endowed with the factor space topology is of the type $\BDiff M_{\alpha}$.  It will be our model of $\BDiff M_{\alpha}$. 

\subsection{A model for maps (\ref{eq:1.1})}\label{ss:5.3} The correspondence $f\mapsto f'$ described in (iv) determines the correspondence $[f]\mapsto [f']$ and therefore determines a specific model for the map (\ref{eq:1.1}) in the case where $M_{\beta}$ is obtained from $M_{\alpha}$ by attaching handles along non-empty boundary components $M_{\alpha}$. 

Next we need a specific model for the map (\ref{eq:1.1}) in the case where the manifold $M_{\beta}$ is obtained from $M_{\alpha}$ by attaching several handles along non-empty components of $\partial M_{\alpha}$ and several handles along empty components of $M_{\alpha}$. To simplify the description of the map, let us assume that $M_{\beta}$ is obtained from $M_{\alpha}$ by attaching $l$ handles along $l$ ordered empty components of $M_{\alpha}$. In fact, let us assume that $M_{\alpha}$ is closed, i.e., $M_{\alpha}$ has no non-empty boundary components. In the general case the definition of the model for maps (\ref{eq:1.1}) is similar. 

\begin{remark} In the general case, however, we essentially use our assumption that if $\partial_kM$ is an empty component of $(M, w)$, then all components $\partial_lM$ with $l<k$ are also empty. 
\end{remark}

%Suppose $M_{\alpha}$ has $k$ ordered empty boundary components $S^{-1}\times S^{d}$, denoted by $\partial_1M_{\alpha}, %\partial_2M_{\alpha}, \dots, \partial_kM_{\alpha}$.  Let $M_{\beta}$ be the manifold obtained from $M_{\alpha}$ by %attaching $l$ handles $D^0\times S^d$ along empty components indexed $\alpha_1, ..., \alpha_l$. Let us note that the %ordered empty boundary components of $M_{\beta}$ correspond to boundary components of $M_{\alpha}$ in a way so that %the component $\partial_iM_{\beta}$ for $1\le i\le k-l$ corresponds to $\partial_jM_{\alpha}$ where %$\alpha_{j-i}<j<\alpha_{j-i+1}$.   

Let $f_{\alpha}: M_{\alpha}\to \R^{\infty}$ be a map representing a point $[f_{\alpha}]$ in our model of $\BDiff M_{\alpha}$. We will construct a map $f_{\beta}: M_{\beta}\to \R^{\infty}$ representing the image $[f_{\beta}]$ of $[f_{\alpha}]$ in our model of $\BDiff M_{\beta}$. 

Let $r(f_{\alpha})$ denote the minimal radius of a closed disc $D_{r(f_{\alpha})}$ in $\R^{\infty}$ centered at $0\in\R^{\infty}$ for which there is a decomposition of the manifold $M_{\alpha}$ into a disjoint union of two submanifolds $M_{\alpha}^1$ and $M_{\alpha}^2$ such that  
\begin{itemize}
\item the disc $D_{r(f_{\alpha})}$ is the minimal round disc centered at $0\in \R^{\infty}$ that contains the image $f_{\alpha}(M^1_{\alpha})$, and 
\item each path component $S$ of $M^2_{\alpha}$ is diffeomorphic to a sphere $S^d$ and the image $f_{\alpha}(S)$ is given by the unit sphere 
\[
    D^{0}\times S^{d} \subset \R^0\times \R^{d+1}\times \{10j+r(f_{\alpha}) \}\times \{0\}\times \{0\}  
\]
for some positive integer $j$.
\end{itemize}

The map $f_{\beta}$ is defined to be the disjoint union of $f$ and the inclusion of $l$ almost standard tori
\[
    D^{0}\times S^{d} \subset \R^0\times \R^{d+1}\times \{10(i+u(i))+r(f_{\alpha}) \}\times \R\times \{0\},  
\]
where $i\in\{1,2,..., k\}$ is the ordering index of an attached handle and $u(i)$ is the number of spheres $f_{\alpha}(S)$ with $S\in M^2_{\alpha}$ whose centers are in the closed disc in $\R^{\infty}$ centered at zero of radius $10i+r(f)$. 

We need to verify that the maps 
\[      
      \BDiff M_{\alpha}\longrightarrow \BDiff M_{\beta}
\]
\[
        [f_{\alpha}]\mapsto [f_{\beta}]
\]
form a strictly commutative diagram, i.e., any two sequences of maps in the diagram
\[  
     \BDiff M_{\alpha}=\BDiff M_0 \longrightarrow \BDiff M_1 \longrightarrow \cdots \longrightarrow \BDiff M_p=\BDiff M_{\beta},
\]
and
\[  
     \BDiff M_{\alpha}=\BDiff M'_0 \longrightarrow \BDiff M'_1 \longrightarrow \cdots \longrightarrow \BDiff M'_q=\BDiff M_{\beta},
\]
where $p$ and $q$ are two positive (possibly different) integers, take each point $[f_{\alpha}]$ onto the same point $[f_{\beta}]$. To this end let us observe that if the submanifold $M^2_{\alpha}$ is not empty, then $[f_{\alpha}]$ belongs to the image of an embedding
\begin{equation}\label{eqn:1}
   \BDiff M_{-1}\longrightarrow \BDiff M_{\alpha}
\end{equation}
where $M_{-1}$ is diffeomorphic to $M^1_{\alpha}$ and therefore, by precomposing the two sequences of maps with the embedding (\ref{eqn:1}), we may assume that $M^2_{\alpha}$ is empty. In this case, it is easy to see that the resulting map $f_{\beta}$ maps the sphere $S$ attached along the $j$-th component onto 
\[
    D^{0}\times S^{d} \subset \R^0\times \R^{d+1}\times \{10j+r(f_{\alpha}) \}\times \{0\}\times \{0\}  
\]
and therefore the resulting map does not depend on the choice of one of the two sequences of maps.

\section{Definition of the homotopy colimit $\mathfrak{H}(d)$}\label{s:6}

Next we describe the homotopy colimit diagram of $\mathfrak{H}(d)$ which is a covariant functor $\mathcal{D}$ from a small category $\I$ to the category of topological spaces. The objects in $\I$ are in bijective correspondence with equivalence classes $[M, w]$ of pairs of an oriented compact manifold $M$, one for each diffeomorphism type, with labeled boundary components, each of the form $S^{i-1}\times S^{d-i}$, and parametrizations $w$ of $\partial M$. 

If $\partial M$ consists of one component $S^{i-1}\times S^{d-i}$, then in $\I$ there are precisely two arrows from the object corresponding to $[M, w]$, namely, an arrow corresponding to attaching a handle $S^{i-1}\times D^{d-i+1}$ and  an arrow corresponding to attaching a handle $D^{i}\times S^{d-i}$, both attachments are by means of the parametrization $w$. In particular, if $\partial M$ is $S^{-1}\times S^d$, then in $\I$ there are two arrows from the object of $[M, w]$; namely, one corresponding to forgetting the empty boundary component $S^{-1}\times S^d$, and one corresponding to taking the disjoint union of $M$ with the sphere $S^d$.  

In general, if $M$ is bounded by $n$ components of the form $S^{i-1}\times S^{d-i}$, then in the category $\I$ there are precisely 
\[
3^n-1=\sum_{k=1}^{n} 2^k \frac{n!}{(n-k)!k!}
\] 
arrows from the vertex corresponding to $[M, w]$; each arrow corresponds to an attachment of handles along some of the components of $\partial M$. 

More precisely, let $[M, w]$ be a pair corresponding to an object in $\I$. Suppose that it is represented by a compact oriented manifold $M$ with $n$ labeled boundary components $\partial_1M, ..., \partial_nM$, each of the form $S^{i-1}\times S^{d-i}$; and a parametrization $w$ of $\partial M$. Then each regular arrow from the object corresponding to $[M, w]$ is labeled by 
\begin{itemize}
\item a choice of $k>0$ components of $\partial M$, say $\partial_{j_1}M, ..., \partial_{j_k}M$, and 
\item one of the two choices of handles $D^{i}\times S^{d-i}$ or $S^{i-1}\times D^{d-i+1}$ for each chosen boundary component of $\partial M$ of the type $S^{i-1}\times S^{d-i}$.    \end{itemize}
Such an arrow points to an object corresponding to a pair $[M', w']$, which can be described as follows. The manifold $M'$ is obtained from $M$ by attaching the chosen handles along $\partial_{j_1}M, ..., \partial_{j_k}M$ by means of $w$. We note that there is a canonical inclusion of $M$ into $M'$. In particular, the path components of $\partial M'$ inherit an order from the order of boundary path components of $M$. We label the path components of $\partial M'$ according to the inherited order. The parametrization $w'$ is defined to be the one given by the composition of the map $\partial M'\to \partial M$ inverse to the canonical inclusion and $w$. Clearly, the equivalence class $[M', w']$ constructed above does not depend on the choices of $M$ and $w$ representing $[M, w]$.

\begin{remark}\label{r:5.1}  There is no canonical choice of $M_{\alpha}$, $w_{\alpha}$ representing $[M_{\alpha}, w_{\alpha}]$ and of $M_{\beta}, w_{\beta}$ representing $[M_{\beta}, w_{\beta}]$. In particular, the existence of an arrow in the category $\I$ from an object of $[M_{\alpha}, w_{\alpha}]$ to an object of $[M_{\beta}, w_{\beta}]$ does not imply the existence of a canonical inclusion of $M_{\alpha}$ to $M_{\beta}$.  
\end{remark}

Next we define the covariant functor $\mathcal{D}$ from the category $\I$ to the category of topological spaces. 

Informally, if $\alpha$ is an object of $\I$ corresponding to a pair $(M_{\alpha}, w_{\alpha})$, then $\mathcal{D}(\alpha)$ is the topological space given by the classifying space $\BDiff M_{\alpha}$ of the group of orientation preserving diffeomorphisms of $M_{\alpha}$ trivial on the boundary. Since an arrow $g: \alpha\to \beta$ in $\I$ corresponds to attaching handles to a manifold $(M_{\alpha}, w_{\alpha})$ along components of the parametrization $w_{\alpha}$, there exists an inclusion of $M_{\alpha}$ into the resulting manifold $M_{\beta}$. The morphism $\mathcal{D}(g)$ is defined to be the induced map
\begin{equation}\label{eq:-1}
     \BDiff M_{\alpha}\longrightarrow \BDiff M_{\beta}. 
\end{equation}
Since the inclusion $M_{\alpha}$ to $M_{\beta}$ is not unique (see Remark~\ref{r:5.1}) some care is necessary in order to get a commutative diagram of spaces (see Remark~\ref{re:5.1}). Namely we choose the models of classifying spaces $\BDiff M_{\alpha}$ and maps (\ref{eq:-1}) as defined in section~\ref{s:model}. Then the diagram of maps (\ref{eq:-1}) is commutative and therefore the covariant functor $\mathcal{D}$ is well defined.

Finally, the space
\[
    \mathfrak{H}(d)\colon= \mathop{\mathrm{hocolim}}_{\alpha\in \mathfrak{I}}\ \BDiff M_{\alpha}
\]
is called the {\it homotopy colimit $\mathfrak{H}(d)$ of classifying spaces} of diffeomorphism groups of manifolds of dimension $d$.

\begin{remark}\label{rem:1} Let $M$ be an oriented manifold with boundary of the form $S^{i-1}\times S^{d-i}$. Given a parametrization $w$ of $\partial M$, the pair $(M, w)$ represents a vertex in the homotopy colimit diagram of $\mathfrak{H}(d)$. There are two arrows from $(M, w)$, namely, the one to a vertex of the manifold $M_{0}$ obtained from $M$ by attaching a handle $D^{i}\times S^{d-i}$ along $w$, and the one to a vertex of the manifold $M_1$ obtained from $M$ by attaching a handle $S^{i-1}\times D^{d-i+1}$. In other words, in the homotopy colimit diagram of $\mathfrak{H}(d)$ there is a subdiagram
\[
\xymatrix
{
& \BDiff M_0 &       &\BDiff M_1     \\     
&            & \BDiff M \ar[ul] \ar[ur] & 
}
\]
If $d$ is even, then the target spaces $\BDiff M_0$ and $\BDiff M_1$ in the above diagram are always different. However if $d$ is odd, the target spaces may coincide. In particular the term $\mathop{\mathrm{colim}}^1 H_0$ in the homology Bousfield-Kan spectral sequence contributes a cycle in $\mathfrak{H}(d)$ of degree $1$.   
\end{remark} 
  
%\begin{example}\label{ee:1}
%If $M$ is a manifold diffeomorphic to $S^1\times D^2$, then there is a parametrization %$w$ of $\partial M$ such that both manifolds $M_0$ and $M_1$ in Remark~\ref{rem:1} are %diffeomorphic to $S^3$. As has been explained in Remark~\ref{rem:1}, this produces a %non-trivial phantom cycle in $\mathfrak{H}(d)$. 
%\end{example}

\begin{example}\label{ex:6.6} In the case where $d=0$, an oriented compact manifold of dimension $d$ is a disjoint union of finitely many points each of which is labeled by a sign, ``+" or ``--". In this case a manifold may have only empty boundary components. In particular we may omit mentioning parametrizations of boundary components. However, the order of boundary components is still essential. 

Let $M=M^b_{p,n}$ be an oriented manifold that consists of $p$ positive points, $n$ negative points, and $b$ empty boundary components $S^{-1}\times S^0$. An orientation preserving diffeomorphism of $M$ is an element of the group $\Sigma_{p,n}$ of pairs of permutations of $p$ and $n$ points. Our model for $\BDiff M=\mathrm{B}\Sigma_{p,n}$ is the standard one; it is defined to be the space of subsets of unordered $p+n$ points in $\R^{\infty}$ such that $p$ points are labeled by ``+" and $n$ points are labeled by ``--". The integer $p-n$ is called the \emph{charge} of the manifold $M$. 

The category $\I$ has one object for each diffeomorphism type of oriented manifolds $M$ of $p$ positive points, $n$ negative points and $b$ empty boundary components, i.e., one object for each triple of non-negative integers $b,p$ and $n$. Attaching $s$ non-empty handles to an oriented manifold $M$ of dimension $0$ results in taking a disjoint union of $M$ and $s$ copies of a pair of a positive point and a negative point. It follows that the category $\I$ splits into the disjoint union of subcategories, one for each integer $\mu$, where the subcategory indexed by an integer $\mu$ is a full subcategory of $\I$ that consists of all objects of manifolds of charge $\mu$.  

Let us fix an integer $\mu$, and let $\J$ denote the full subcategory of $\I$ indexed by $\mu$. In the category $\J$, each arrow is a composition of \emph{elementary} arrows that are arrows from an object with $b$ boundary components to an object with $b-1$ boundary components. 
 Furthermore, from an object $M^b_{p,n}$, there are precisely $b$ elementary arrows in $\J$ to an object of $M^{b-1}_{p+1, n+1}$, and $b$ elementary arrows to an object of $M^{b-1}_{p,n}$, and there are no other elementary arrows from the object of $M_{p,n}$. 
 
Let $\J(b)$ denote the full subcategory of $\J$ that consists of objects of manifolds with $\le b$ boundary components. Then, clearly,  
\[
   \mathop{\mathrm{hocolim}}_{\alpha\in \J(0)} \BDiff M_{\alpha} =\bigsqcup_{p-n=\mu} {B\Sigma}_{p,n}. 
\] 
It is also easy to observe that
\begin{equation}\label{equ:6.1}
   \mathop{\mathrm{hocolim}}_{\alpha\in \J(1)} \BDiff M_{\alpha} =B\Sigma_{\mu+\infty,\infty}\simeq B\Sigma_{\infty}\times B\Sigma_{\infty},
\end{equation}
where $B\Sigma_{\mu+\infty, \infty}$ is the homotopy colimit of the diagram
\[
   B\Sigma_{\mu, 0}\longrightarrow B\Sigma_{\mu+1,1}\longrightarrow B\Sigma_{\mu+2,2} \longrightarrow \cdots
\]
if $\mu\ge 0$, and 
\[
   B\Sigma_{0, \mu}\longrightarrow B\Sigma_{1, \mu+1}\longrightarrow B\Sigma_{2, \mu+2} \longrightarrow \cdots
\]
if $\mu<0$. The construction of 
\[ 
  \mathop{\mathrm{hocolim}}_{\alpha\in \J(\infty)} \BDiff M_{\alpha}
\]
resembles the Quillen ``+'' construction applied to the space (\ref{equ:6.1}).

\end{example}

\section{Elementary bordisms}\label{se:7}

\subsection{Standard handle of index $i$.}
Let $d$ be an arbitrary fixed non-negative integer. For $i\in\{0,1, \dots, \left\lfloor \frac{d+1}{2}\right\rfloor\}$, let $m_i: \R^{d+1}\to \R$ be the standard Morse function of index $i$ given by
\[
  (x_1, ..., x_{d+1})\mapsto -x_1^2-\cdots -x_i^2+x_{i+1}^2+\cdots +x_{d+1}^2.
\]
Then there is an open neighborhood $U'$ of the union of coordinate planes
\[
   x_1^2+\cdots + x_i^2 = 0 \qquad \mathrm{and} \qquad x_{i+1}^2+\cdots + x_{d+1}^2=0
\]
homeomorphic to a disc such that the closure $U$ of the set
\[
   U'\cap \{\, x\in \R^{d+1}\, |\, -1\le m_i(x)\le 1\, \}
\]
is a {\it handle of index $i$}  diffeomorphic to a closed disc with corners; the set of corner points is empty for $i=0$. The boundary $\partial U$ of $U$ is the union $\partial U=U_h\cup U_v$ of two components  
\[
   U_h=\overline{U'\cap \{\,m_i(x)=\pm 1 \,\}} \quad \mathrm{and}\quad U_v\cong  [-1,1]\times S^{i-1}\times S^{d-i};
\] 
the component $U_v$ is empty for $i=0$. Without loss of generality we may assume that $U'$ is chosen so that $U_v$ and $U_h$ are smooth submanifolds of $\R^{d+1}$. 

We fix the latter diffeomorphism parametrizing $U_v$ and consider it to be a part of the structure of the handle $U$. We will assume that the diffeomorphism parametrizing $U_v$ is orientation preserving. 

The function $m_i|U: U\to [-1,1]$ is said to be the {\it standard model} of a Morse function on a handle of index $i$ if 
$U_v$ is an integral submanifold of $\R^{d+1}$ with respect to the gradient vector field of $m_i$, and the projection of $U_v$ onto its first factor coincides with $m_i|U_v$.

\subsection{Kernel of an elementary bordism}

Let $W$ be an oriented bordism between oriented compact manifolds $-M_0$ and $M_1$ of dimension $d$, i.e., $W$ is an oriented compact manifold of dimension $d+1$ with boundary given by the disjoint union of oriented manifolds $ -M_0\sqcup  M_1$. It is known that $W$ admits a Morse function $g: W\to [-1,1]$ with 
\[
    g^{-1}(-1)=M_0 \qquad \mathrm{and} \qquad g^{-1}(1)=M_1.
\]
We say that a bordism $(W, g)$ is {\it elementary} if $g$ has a unique critical point $p$ of index $i\in \{0, ... ,\left\lfloor \frac{d+1}{2}\right\rfloor\}$ and $g(p)=0$. We will also say that $(W, g)$ is an elementary bordism of index $i$. We note that in the case where $d$ is even, if $(W, g)$ is elementary, then $(W, -g)$ can not be elementary. A {\it kernel} of an elementary bordism $(W, g)$ is an orientation preserving embedding 
\[
     \psi: U\longrightarrow W
\]
of a standard handle of index $i$ such that the composition $g\circ\psi$ coincides with the standard model $m_i$ of a Morse function on $U$. 

An elementary bordism $(W_1, g_1)$ with kernel $\psi_1$ is {\it equivalent} to an elementary bordism $(W_2, g_2)$ with kernel $\psi_2$ if there is an orientation preserving diffeomorphism $\alpha: W_1\to W_2$ such that 
\[
      g_2=g_1\circ \alpha^{-1} \qquad \mathrm{and} \qquad \psi_2=\alpha\circ \psi_1.
\]
In this case the diffeomorphism $\alpha$ is said to be an equivalence of $(W_1, g_1, \psi_1)$ and $(W_2, g_2, \psi_2)$. We note that the boundary components of $W_1$ and $W_2$ are ordered and an equivalence preserves the order. The equivalence class of $(W, g, \psi)$ is denoted by $[W, g, \psi]$.

\subsection{The elementary bordism associated to a parametrization}

Let $M$ be a compact oriented manifold bounded by a manifold diffeomorphic to $S^{i-1}\times S^{d-i}$ for some  $i\in\{0,1, \dots, \left\lfloor \frac{d+1}{2}\right\rfloor\}$. Then each parametrization $w$ of $\partial M$ determines an oriented elementary bordism $(W, g)$ and a kernel $\psi$. 

Indeed, such a bordism $W$ is given by a smooth orientable compact manifold (without corners) obtained by attaching the standard handle $U$ of index $i$  to $[-1,1]\times M$ along the orientation reversing map
\[
    \id_{[-1,1]}\times w^{-1}\colon U_v\cong [-1, 1]\times S^{i-1}\times S^{d-i}\longrightarrow [-1,1]\times \partial M.
\]

We note that by the definition of the standard model $m_i$ of a Morse function, the projection 
\[
    \pi: [-1,1]\times \partial M\longrightarrow [-1,1]
\]
onto the first factor coincides with the composition of $\id_{[-1,1]}\times w$ and $m_i|U_v$. Consequently, in view of the parametrization $w$, there is a canonical extension of $\pi$ to a Morse function $g$ on $W$. Clearly, the inclusion $U\to W$ is a kernel of $(W, g)$. 

It remains to describe the orientations of $W$, $M_0$ and $M_1$. On $W$ we choose an orientation so that at each point of 
\[
M= \{0\}\times M\subset [-1,1]\times M
\] in $W$ a vector frame $\{v_1, ..., v_{d+1}\}$ in $TW$ is positive, where $v_1$ maps by $dg$ to the vector of the positive direction in $[-1,1]$ and the frame $\{v_2, ..., v_{d+1}\}$ is a frame in $TM\subset TW$ corresponding to the positive orientation of $TM$. Next we choose orientations on $M_0$ and $M_1$ so that $W$ is an oriented bordism between $-M_0$ and $M_1$. We warn the reader that our choice of orientation of $M_0$ is not consistent with the canonical orientation\footnote{Let us recall that if $N$ is an oriented manifold with boundary, then the boundary $\partial N$ inherits an orientation: a frame $\{v_2, ..., v_{d+1}\}$ at $p\in \partial N$ is an oriented basis of $T_p\partial N$ if and only if the frame $\{v_1, ..., v_{d+1}\}$ is an oriented basis of $T_pN$, where $v_1$ is a vector outward normal to $N$.} of the boundary $\partial W$. 

In fact each equivalence class $[M, w]$ of a compact oriented manifold $M$ with $\partial M$ of the form $S^{i-1}\times S^{d-i}$ and a parametrization $w$ of $\partial M$ gives rise to a well-defined equivalence class $[W, g, \psi]$ of an elementary bordism $(W, g)$ with a kernel $\psi$. Indeed, an equivalence  $\alpha: (M_1, w_1)\to (M_2, w_2)$ of pairs representing the same class $[M, w]$ determines a diffeomorphism 
\[
    \id_{[-1,1]}\times \alpha: [-1,1]\times M_1\longrightarrow [-1,1]\times M_2
\]  
which canonically extends to a diffeomorphism $W_1\to W_2$ of elementary bordisms associated with $(M_1, w_1)$ and $(M_2, w_2)$ respectively. Furthermore, the obtained diffeomorphism $W_1\to W_2$ is an equivalence of elementary bordisms.

\begin{remark}\label{rem:7.1} For $i=1,2$, let $(W_i, g_i, \psi_i)$ be an elementary bordism obtained from a pair $(M_i, w_i)$. As we have seen, an equivalence of pairs $(M_i, w_i)$ produces an equivalence of triples $(W_i, g_i, \psi_i)$. On the other hand most of the equivalences of triples $(W_i, g_i, \psi_i)$ are not produced by equivalences of pairs $(M_i, w_i)$.  
\end{remark}

\subsection{The parametrization associated to an elementary bordism}\label{ss:3.3}
Let $(W, g)$ be an oriented elementary bordism between oriented closed manifolds $M_0$ and $-M_1$ and $\psi$ be a kernel of $(W, g)$. Then the triple $(W, g, \psi)$ gives rise to an oriented manifold $M$ with boundary of the form $S^{i-1}\times S^{d-i}$ and an equivalence class $[w]$ of parametrizations of $\partial M$. Indeed, such a 
manifold $M$ is defined to be the submanifold
\[
    g^{-1}(0)\ \cap\ \overline{W\setminus \psi(U)}
\]
of $W$ with orientation that agrees with the canonical orientation of the boundary of the manifold with corners $g^{-1}([-1,0])\cap (W\setminus \psi(U))$;  
 while the equivalence class $[w]$ is represented by the parametrization $w$  given by the identification of $\partial M$ with $\{0\}\times S^{i-1}\times S^{d-i}$ by means of
\[
     \psi|U_v\colon U_v\cong [-1,1]\times S^{i-1}\times S^{d-i}\longrightarrow [-1,1]\times \partial M.
\]

In contrast to the statement of Remark~\ref{rem:7.1}, each equivalence of elementary bordisms with kernels descends to an equivalence of the associated pairs of manifolds $M$ and parametrizations of $\partial M$. It is easily verified that equivalence classes of triples $(W, g, \psi)$ are in bijective correspondence with equivalence classes of pairs $(M, w)$. We will prove a more general statement below (see Lemma~\ref{l:4.1}).   

\begin{remark} To bring our discussion in perspective, let us observe that in the diagram of the homotopy colimit $\mathfrak{H}(d)$ there is precisely one subdiagram
\[
\xymatrix
{
& \BDiff M_0 &       &\BDiff M_1     \\     
&            & \BDiff M \ar[ul] \ar[ur] & 
}
\]
for each equivalence class of triples $(W, g, \psi)$ of an oriented elementary bordism between manifolds $M_0$ and $M_1$, and a kernel $\psi$.  
\end{remark}

Now we can present an example for a claim in Remark~\ref{rem:1}. 

\begin{example}\label{e:3} Let $f$ be a Morse function on $\C P^2$ with precisely three critical points at which the function $f$ assumes values $-2, 0$ and $2$ respectively. Then the manifold 
\[
     W\colon = \C P^2 \setminus \{\, x\in \C P^2\, |\,  -1\le f(x) \le 1\, \}
\]
determines an elementary bordism $(W, f|W)$ between two copies of $S^3$. There are two non-equivalent choices of a kernel of $(W, f|W)$. Let $\psi$ denote one of them. Then, the equivalence class of the triple $(W, f|W, \psi)$ determines a subdiagram 
\[
     \BDiff S^{3} \longleftarrow \BDiff M\longrightarrow \BDiff S^3
\]
of the homotopy colimit diagram of $\mathfrak{H}(d)$. 
\end{example}

\section{Elementary $k$-bordisms}
\subsection{Standard $(k,s,i)$-handle.} 
\begin{definition} A point $x\in M$ of a smooth map $f: M\to N$, with $\dim M\ge \dim N$, is said to be {\it regular} if the kernel rank of $f$ at $x$ is precisely $\dim M - \dim N$. Otherwise, i.e., if the kernel rank of $f$ at $x$ is greater than $\dim M-\dim N$, the map $f$ is said to be \emph{singular} at the point $x$. In this case we will also say that $x$ is a \emph{singular} point of $f$.  

Given a smooth map $f:M\to N$, we say that a fiber $f^{-1}(y)$ over a point $y\in N$ is \emph{regular} if each point $x\in f^{-1}(y)$ is regular. Otherwise, the fiber $f^{-1}(y)$ is said to be \emph{singular}.  
\end{definition}

\begin{definition}
We say that a point $x\in M$ is a fold singular point of a smooth map $f: M\to N$ if there are coordinate neighborhoods of $x$ in $M$ and $f(x)$ in $N$ with respect to which $f$ is of the form
\begin{equation}\label{eq:2.1}
     (x_1, ..., x_{m+d}) \mapsto (x_1, ..., x_{m}, \pm x^2_{m+1} \pm \cdots \pm x^2_{m+d}). 
\end{equation} 
\end{definition}

It immediately follows that if $f: M\to N$ is a fold map, then the singular set $S=S(f)$ of $f$ is a submanifold of $M$ of dimension $\dim N-1$. Furthermore, the restriction of $f$ to $S$ is an immersion. If the manifold $M$ is of dimension $m+d$, then in a neighborhood of each point $x\in S$ there are coordinates with respect to which $f$ has the form (\ref{eq:2.1}) and the number $i$ of signs $``-"$ in the quadratic expression is at most $\left\lfloor \frac{d+1}{2}\right\rfloor$. The number $i$ is called the \emph{index} of the fold singularity of $f$ at $x$. We will also say that the singular point $x$ has index $i$. We note that any two points of the same path component of $S$ have the same index.

We recall that if $X$ is an oriented manifold, then $-X$ stands for an oriented manifold whose underlying space is the same as that of $X$, but the orientations of $X$ and $-X$ do not agree. Let $m(k,s,i)$ denote the map given by 
\[
     (-1)^{k-s}\R^{d+k}\longrightarrow \R^k,
\]
\[
    (y_1,...,y_{k-1}, x_1, ..., x_{d+1})\mapsto (y_1, ..., y_{s-1}, m_i(x_1, ..., x_{d+1}), y_s, ..., y_{k-1})
\]
for 
\[
   k=1, 2, ..., \qquad s=1,2, \dots, k, \quad \textrm{and} \quad i=0,1, \dots, \left\lfloor \frac{d+1}{2}\right\rfloor.
\]
Then the restriction of $m(k,s,i)$ to $(-1)^{k-s}D^{k-1}\times U$, where $U$ is the standard handle of index $i$, is called the standard model on a $(k,s,i)$-handle $(-1)^{k-s}D^{k-1}\times U$. We note that the orientation of $(-1)^{k-s}D^{k-1}\times U$ coincides with the product of standard orientations on $D^{k-1}$ and $U$ only if $k-s$ is even. 

We note that each map $m(k,s,i)$ is a fold map. The singular set of $m(k,s,i)$ is given by $x_1=\cdots =x_{d+1}=0$ and the image of the singular set is given by the set of points of the form $(y_1, ..., y_{s-1}, 0, y_s, ..., y_{k-1})$. Each singular point of $m(k,s,i)$ has index $i$.

\subsection{Parametrization of a multicomponent boundary.}

From each vertex of the homotopy colimit diagram of $\mathfrak{H}(d)$ there are two elementary arrows. These are encoded by an elementary bordism with a kernel. To encode arrows in the diagram of $\mathfrak{H}(d)$ given by compositions of $k$ elementary arrows, we need a notion of an elementary $k$-bordism, which generalizes the notion of an elementary bordism.

Let $D^k$ denote a standard disc of dimension $k$. 
We define an {\it elementary $k$-bordism $(W, g)$} to be a smooth proper map $g: W\to D^k$ of an oriented manifold $W$ with boundary $\partial W=g^{-1}(\partial D^k)$ such that the set of singular values of $g$ coincides with the union of coordinate hyperplanes in $D^k$ and each 
singular point of $g$ is a fold.   

A {kernel} of an elementary $k$-bordism $(W, g)$ is defined to be an orientation preserving embedding $\psi=\sqcup_{s=1}^k \psi_s$ into $W$ with components
\[
     \psi_s\colon (-1)^{k-s}D^{k-1}\times U_{s} \longrightarrow W, 
\]
where $U_s$ is a standard handle of index $i$ for some $i\in\{1,...,\left\lfloor \frac{d+1}{2}\right\rfloor\}$, such that for each $s$ the composition 
\[
      (-1)^{k-s}D^{k-1}\times  U_s\stackrel{\psi_s}\longrightarrow W \stackrel{g}\longrightarrow D^k,   
\]
coincides with the standard model on a $(k,s,i)$-handle. 

We say that an elementary $k$-bordism $(W_1, g_1)$ with kernel $\psi_1$ is {\it equivalent} to an elementary $k$-bordism $(W_2, g_2)$ with kernel $\psi_2$ if there is an orientation preserving diffeomorphism $\alpha: W_1\to W_2$ such that
\[
      g_2=g_1\circ \alpha^{-1} \qquad 
      \mathrm{and}\qquad \psi_2=\alpha\circ \psi_1.
\] 

%We say that an elementary $k$-bordism $(W_1, g_1)$ with kernel $\phi$ is {\it %equivalent} to an elementary $k$-bordism $(W_2, g_2)$ with kernel $\psi$ if %there is a (not necessarily orientation preserving) diffeomorphism $\alpha: %W_1\to W_2$ and a diffeomorphism 
%\[
%   \beta: D^k\longrightarrow D^k
%\]
%permuting the (oriented) factors by means of a permutation $\sigma\in %\Sigma_k$ such that
%\[
%      g_2=\beta\circ g_1\circ \alpha^{-1}
%\] 
%and for each $i=1,..., k$ and $j=\sigma(i)$ the kernel $\psi_j$ is given by %the composition 
%\[
%    D^{j-1}\times U_j \times D^{k-j} \longrightarrow D^{i-1}\times U_i \times %D^{k-i} \stackrel{\phi_i}\longrightarrow W \stackrel{\alpha}\longrightarrow %W, 
%\]
%where the first map permutes the factors by means of $\sigma^{-1}$.

\begin{lemma}\label{l:4.1} Equivalence classes of pairs $(M, w)$ of an oriented compact manifold $M$ with $k$ boundary components of the form $S^{i-1}\times S^{d-i}$ and a parametrization $w$ of $\partial M$ are in bijective correspondence with equivalence classes of triples $(W, g, \psi)$ of an elementary $k$-bordism $(W, g)$ and a kernel $\psi$. 
\end{lemma}
\begin{proof}  Let $M$ be an oriented compact manifold with $k$ boundary components of the form $S^{i-1}\times S^{d-i}$ and $w$ a parametrization of $\partial M$ which has $k$ labeled components $w_1,..., w_k$.
We define an oriented manifold 
\[
       W:= D^k\times M\bigcup\ \mathop{\cup}_{s=1}^{k} (-1)^{k-s}D^{k-1}\times U_s,
\]
by attaching $k$ standard $(k,s,i)$-handles $U_s$ to $D^k\times M$. The $s$-th handle is attached to the $s$-th boundary path component $D^k\times \partial_sM$ by means of an orientation reversing diffeomorphism 
\begin{equation}\label{eq:4.1}
    D^{k-1}\times U_{v,s}\longrightarrow D^{k}\times \partial_sM
\end{equation}
where $U_{v,s}$ is the vertical part of $\partial U_s$, given by the composition of 
\begin{itemize}
\item the orientation preserving diffeomorphism 
\[
 (-1)^{k-s}D^{k-1}\times U_{v,s}\cong [(-1)^{k-s}D^{s-1}\times D^{k-s}]\times [D^1\times S^{i-1}\times S^{d-i}],
\] 
that is the product of $\id_{D^{k-1}}$ and the fixed orientation preserving parametrization of $U_{v,s}$, 
\item the orientation reversing diffeomorphism from
\[
(-1)^{(k-s)}D^{s-1}\times D^{k-s}\times D^1\times S^{i-1}\times S^{d-i}
\]
to
\[
(-1)^{(k-s)}D^{s-1}\times D^{k-s}\times D^1\times \partial_s M, 
\]
given by $\id_{D^{s-1}}\times \id_{D^{k-s}}\times \id_{D^1}\times w_s^{-1}$,
\item the orientation preserving diffeomorphism 
\[
(-1)^{(k-s)}D^{s-1}\times D^{k-s}\times D^1\times \partial_s M \longrightarrow 
D^{s-1}\times D^1\times D^{k-s}\times 
\partial_s M
\]
exchanging factors,
\item and the identity  
\[
D^{s-1}\times D^1\times D^{k-s}\times 
\partial_s M\longrightarrow D^k\times \partial_s M.
\]
 \end{itemize}

As in the case of an elementary  bordism, the projection $D^k\times \partial_s M\to D^k$ onto the first factor agrees with the restriction of the standard model
\[
     m(k,s,i): (-1)^{k-s}D^{k-1}\times U_s\longrightarrow D^k
\]
to $(-1)^{k-s}D^{k-1}\times U_{v,s}$, i.e., there is a commutative diagram \[
\xymatrix
{
&  (-1)^{k-s}D^{k-1}\times U_{v,s} \ar[dr] \ar[rr]^{(\ref{eq:4.1})}&       &  D^k\times \partial_sM \ar[dl]     \\     
&            & D^k.   & 
}
\]

 Consequently, the projection $D^k\times M\to D^k$ onto the first factor extends to an elementary $k$-bordism $(W, g)$ with a kernel $\psi$ embedding the union $\cup U_s$ of standard handles into $W$. 

Suppose now that $(W, g)$ is an elementary bordism and $\psi$ is a kernel of $(W, g)$. Then
\[
     M\colon = g^{-1}(0) \ \cap \overline{W\setminus \cup \psi(U_s)}. 
\]
is an oriented compact manifold with $k$ boundary components of the form $S^{i-1}\times S^{d-i}$. Each component of $\partial M$ is identified by means of $\psi$ with $\{0\}\times U_{h,s}$ for some $s$. In particular, there is a canonical parametrization $w$ of $\partial M$. The index $s$ determines an order on the boundary components of $\partial M$. Thus each triple $(W, g, \psi)$ leads to a pair $(M, w)$. 

\end{proof}

\section{Marked fold maps}\label{se:9}

\subsection{Motivation}
Stable characteristic classes of smooth manifold bundles turn out to be related to fold maps equipped with an additional structure, which we call a \emph{marking}. The notion of a marked fold map is a generalization of the notion of a $k$-bordism with a kernel. A marked fold map is a fold map of smooth manifolds with two compatible markings, a global one and a local one. These are motivated by the following observation. 

To begin with let us observe that any fold map near a singular fiber looks like a $k$-bordism. More precisely, let $f: M\to N$ be a fold map whose restriction to the set of singular points is an immersion with normal crossings. Let $\Sigma_k$ denote the set of points $x$ in $N$ such that the fiber of $f$ over $x$ contains precisely $k$ singular points. It immediately follows that if $x\in \Sigma_k$, then the map $f$ restricted to the inverse image $f^{-1}(D)$ of a disc neighborhood $D$ of $x$ admits a decomposition into a product of an elementary $k$-bordism and a constant map of a disc of dimension $\dim N-k$. 
A local marking fits such decompositions for $x$ ranging over points of $\Sigma_k$ with a fixed $k$, while a global marking fits local markings for strata $\Sigma_k$ with different $k$.

\subsection{Local marking}
Let $f: M\to N$ be a fold map whose restriction to the set of singular points is an immersion with normal crossing. Let $\Sigma_k$ denote the set of points $x$ in $N$ for which the fiber of $f$ over $x$ is a submanifold of $M$ with precisely $k$ singular points. Then each $\Sigma_k$ is a submanifold of $N$ of codimension $k$. Suppose that for each $k>0$ and each path component $\Sigma$ of $\Sigma_k$, there is a fixed smooth orientation preserving diffeomorphism 
\[
i_N: D^k\times\Sigma\hookrightarrow N\quad \mathrm{extending}\quad \{0\}\times \Sigma\hookrightarrow N
\] 
and an elementary $k$-bordism $\pi: W\to D^k$ with a kernel $\psi=\{\psi_s\}$ together with a fixed smooth orientation preserving diffeomorphism $i_M: W\times \Sigma \hookrightarrow M$ for which the diagram 
\[
   \begin{CD}
    W\times \Sigma @>i_M>> M \\
    @V\pi\times\id_{\Sigma}VV @Vf VV \\
    D^k\times \Sigma @>i_N>> N
   \end{CD}
\]
commutes. Then we say that the set of pairs of maps $i_M, i_N$, one pair for each path component $\Sigma$, is a \emph{local marking} on the fold map $f$. We observe that the image of $i_N$ is path connected, while the image of $i_M$ has precisely $k$ path components. 

\begin{remark} A local marking on $f$ encodes an order of singular points of $f$ in each fiber. 
\end{remark}

\begin{remark} Not every fold map admits a local marking since an existence of a local marking implies the existence of a trivialization of the normal bundle of each path component of $S(f)\cap f^{-1}(\Sigma)$, where $S(f)$ is the set of singular points of $f$ in $M$. Furthermore, if $f$ is a fold map, then there is a well defined kernel bundle over $S(f)$. The existence of a marking implies that the restriction of the kernel bundle to $S(f)\cap f^{-1}(\Sigma)$ is trivial.  
%On the other hand, if $S(f)\cap f^{-1}(\Sigma)$ has precisely $k$ components, %and  the normal bundle of each path component of $S(f)\cap f^{-1}(\Sigma)$ is %trivial, then $f$ admits a local marking. Indeed, to construct a local marking %on $f$, we may first construct embeddings $i_N: D^k\times \Sigma\to N$ so that %the image of coordinate hyperplanes in $D^k\times \{x\}$ maps to singular %values of $f$ for each $x\in \Sigma$. Then, using the Ehresmann theorem, we %may construct embeddings $i_M$ as required.   
%What if the kernel is not trivial.
\end{remark}

\begin{remark}\label{r:6.1} A fold map may contain many essentially different markings. For example, suppose $f: M\to N$ is a fold map of an orientable manifold into a surface. Then the singular set of $f$ is of dimension $1$. Suppose $f$ contains a circle $\Sigma\subset N$ of singular values of index $0$. Let $(i_M, i_N)$ be a marking of $f$ near $\Sigma$. Let $d$ be the dimension of a non-empty regular fiber of $f$. Then for each element of $\tau\in \pi_1(SO_{d+1})$, there is a new marking $(i'_M, i'_N)$ of $f$ near $\Sigma$ such that the difference between the trivializations given by $i'_N$ and $i_N$ is in the homotopy class of $\tau$.    
\end{remark}

\begin{remark} Even if $f$ is a fold map of oriented manifolds such that $\dim M=\dim N$, it may not admit a local marking. Indeed, let $N$ denote the oriented manifold 
\[
    N\colon=   \R^3\times [0,1] /\sim
\]
where $\sim$ identifies a point $(x,y,z; 0)$ with the point $(z,x,y; 1)$. Of course $N$ is diffeomorphic to $\R^3\times S^1$. There is a fold map $f$ from an oriented manifold $M$ to the manifold $N$ such that the set of singular values of $f$ is 
\[
    \{\ (x,y,z; t)\in \R^3\times [0,1]/\sim \ | \ xyz=0\ \}. 
\]  
Clearly $f$ does not admit a local marking. 
\end{remark}

\begin{remark} If $f: M\to N$ is a fold map with a local marking, then for each $k$ and for each path component $\Sigma$ of $\Sigma_k$ the normal bundle of each path component of $S(f)\cap f^{-1}(\Sigma)$ is trivial. On the other hand, it is easy to construct an example of a fold map $f$ with a local marking such that the normal bundle of the singular set of $f$ is not trivial; e.g., consider a fold map $f: M\to \R^2$ of a non-orientable manifold such that for some path component $\Sigma\subset M$ of the singular set of $f$ the image $f(\Sigma)$ is a figure ``$8$".  
\end{remark}

\begin{remark}\label{rem:9.5} Let $f$ be a proper fold map with a non-empty regular fiber of dimension $d$. Then, by Lemma~\ref{l:4.1}, a local marking associates an oriented compact manifold $M$ of dimension $d$ with a parametrized boundary $\partial M$ to each point of the target of $M$. Indeed, if $x$ is a regular value of $f$, then the associated manifold $M$ is $f^{-1}(x)$. If $x$ is a singular value, then $x\in \Sigma$ and there is a fixed elementary $k$-bordism $\pi: W\to D^k$ which is a part of the local marking. By Lemma~\ref{l:4.1}, it corresponds to a manifold of dimension $d$ with a parametrized boundary. 
\end{remark}

\subsection{Global marking}

Let $f:M\to N$ be a fold map of oriented manifolds such that the restriction of $f$ to the set of singular points in $M$ is an immersion with normal crossings. Suppose that for each path component $S\subset M$ of singular points of $f$ of index $i$, there is 
\begin{itemize}
%\item a closed neighborhood of the submanfiold $\Sigma$ in $M$ given by a %submanifold $V\subset M$ with corners,
\item a diffeomorphism $\alpha: U\times S \to M$ extending the inclusion of the singular set $\{0\}\times S$ into $M$, where $U$ is the standard handle of index $i$, and
\item a submersion $\beta: [-1,1]\times S\to N$ 
\end{itemize}
such that the diagram 
\[
\begin{CD}
   U\times S @>>> M \\
   @Vm_i\times \id_{S} VV @Vf VV \\
   [-1,1]\times S @>\beta>> N,
\end{CD}
\]
is commutative. Here $m_i:U\to \R$ is the standard model of a Morse function of index $i$. 
Then the couple $(\alpha, \beta)$ is said to be a {\it kernel} of the fold map $f$. 

\begin{remark}\label{r:9.5} The existence of a kernel implies that the normal bundle of $S$ in $M$ is trivial and the normal bundle of the immersion $f|S: S\to N$ is trivial. Furthermore, a kernel determines a trivialization of the normal bundles of $S$ and $f|S$.  
\end{remark}

\begin{remark} If the manifold $M$ is oriented, then the trivialization of the normal bundle of $S$ determines an orientation on $S$. On the other hand, if $N$ is oriented, then $f(S)$ inherits an orientation since the normal bundle of the immersion $f|S$ is trivial. However, $f|S$ may not preserve the orientation. Indeed, if $f|S$ preserves the orientation, then for the map $\tilde f: M\to (-1)N$ that pointwise coincides with $\tilde f$, the immersion $\tilde f|S$ does not preserve the orientation. 
\end{remark}

\begin{remark} We emphasize that according to our definition, the index $i$ of singular points is at most $\left\lfloor \frac{d+1}{2}\right\rfloor$. In particular, the coorientation of $f(S)$ in $N$ determined by $\beta$ does not really depend on the choice of $\beta$. 
\end{remark}
 Suppose now that for each point $x\in N$ in $f(S)$, the set of singular points of $f$ in $f^{-1}(x)$ is ordered. We say that the order of singular points $S\cap f^{-1}(x)$ is \emph{preserved} as $x$ ranges over singular values of $f$ in $N$ if each point $x\in N$ has a neighborhood $O(x)$ such that for any point $y$ in $O(x)$ the orders of singular points in $S\cap f^{-1}(x)$ and $S\cap f^{-1}(y)$ are compatible. More precisely, each point of $S\cap f^{-1}(y)$ belongs to a unique path component of $S\cap f^{-1}(O(x))$ and therefore corresponds to a unique point in $S\cap f^{-1}(x)$ provided that $O(x)$ is chosen to be small enough. Thus there is a map
 \[
      S\cap f^{-1}(y)\longrightarrow S\cap f^{-1}(x).
 \] 
 We require that this map preserves the order of singular points. 
 
\begin{definition} A \emph{global marking} on a fold map $f: M\to N$ of oriented manifolds is a choice of a kernel $(\alpha, \beta)$ and a choice of an order of singular points in the singular fibers $f^{-1}(x)$ of $f$ such that the order of singular points is preserved as $x$ ranges over singular values of $f$.  
 \end{definition}

\begin{example} An elementary $k$-bordism with a kernel has a canonical global marking. Indeed, the kernel of an elementary $k$-bordism $f$ determines a kernel of $f$ considered as a fold map. Furthermore, for each singular value $x$ of $f$, each singular point in the fiber $f^{-1}(x)$ belongs to a unique $(k,s,i)$-handle $U_s$ and therefore is marked by a positive integer $s$. In particular, the singular points in each singular fiber of $f$ are ordered. It is easy to see that the order of singular points in the fibers $f^{-1}(x)$ is preserved as $x$ ranges over singular values of $f$. 
\end{example}

%is an example of a premarked fold map. On the other hand if $f:M\to N$ is a %premarked fold map, then for any point $x\in N$ there is a coordinate %neighborhood $D$ about $x$ such that 
%\[
%   f|f^{-1}(D)\colon f^{-1}(D)\longrightarrow D
%\]
%is a $k$-bordism with a kernel. 
%\end{example}

\begin{remark} The example in Remark~\ref{r:6.1} shows that a fold map may have many essentially different global markings. 
\end{remark}

\subsection{Definition of a marked fold map}

If a fold map $f$ has a local marking and a global marking, and if these two structures agree in a sense defined below, then we say that $f$ is a marked fold map. 

To describe the compatibility condition, suppose $f$ is a fold map with a local marking and with a global marking. Let 
\[
i_M: W\times \Sigma\to M \qquad \mathrm{and}\qquad i_N: D^k\times \Sigma\to N
\] 
be embeddings of a local marking. The structure of a local marking also includes a projection $g: W\to D^k$ and a kernel $\psi$ of the elementary $k$-bordism $(W, g)$ with components
\[
     \psi_s: (-1)^{k-s} D^{k-1}\times U_s\longrightarrow W. 
\]
We note that the image of $i_M$ of the $s$-th singular set of $g$ is a submanifold of codimension $k-1$ of a path component $\tilde\Sigma$ of fold singular points of $f$. Let 
\[
\alpha: U\times \tilde\Sigma \to M \qquad \mathrm{and}  \qquad \beta: [-1,1]\times \tilde\Sigma\to 	N
\] 
be kernel maps of $f$. Then $U\cong U_s$ and all mentioned maps fit a commutative diagram
\[
   \begin{CD}
     D^{k-1}\times U_s\times \Sigma @>\psi_s\times\id_{\Sigma}>> W\times \Sigma @>i_M>> M @<\alpha<< U\times \tilde\Sigma \\
     @Vm(k,s,i)\times \id_{\Sigma}VV @Vg\times\id_{\Sigma}VV @VfVV @Vm_i\times \id_{\tilde\Sigma} VV \\
     D^k\times \Sigma @>\cong>> D^k\times \Sigma @>i_N>> N @<\beta<< D^1\times\tilde\Sigma,
   \end{CD}
\]
where all horizontal maps are diffeomorphisms onto image. 

We say that the global and local markings of $f$ are \emph{compatible} if for each pair $(i_M, i_N)$ and each pair $(\alpha, \beta)$, there exists a diffeomorphism onto image
\[
   \sigma: D^{k-1}\times \Sigma \longrightarrow \tilde\Sigma
\]
such that the composition $\alpha^{-1}\circ i_M\circ (\psi_s\times \id_{\Sigma})$ of the upper horizontal homomorphisms coincides with
\[
   D^{k-1}\times U_s\times \Sigma \longrightarrow U_s\times (D^{k-1}\times \Sigma) \stackrel{id_{U_s}\times \sigma}\longrightarrow U \times \tilde\Sigma
\] 
and the composition $\beta^{-1}\circ i_N$ of the lower horizontal homomorphisms coincides with 
\[
   D^k\times \Sigma \longrightarrow D^{s-1}\times D^1\times D^{k-s}\times \Sigma \longrightarrow  D^1\times (D^{k-1} \times \Sigma)\stackrel{\id_{D^1}\times\sigma}\longrightarrow  D^1\times\tilde\Sigma.
\]

\begin{definition} A fold map $f$ with a global and local markings is said to be a \emph{marked fold map} if its markings are compatible. 
\end{definition}

\begin{remark} Let $\mathbf{F}(M,N)$  (respectively, $\mathbf{MF}(M,N)$) denote the topological space of fold maps (respectively, marked fold maps) of a manifold $M$ into a manifold $N$. Then there is an inclusion 
\[
   i\colon \mathbf{MF}(M,N)\longrightarrow \mathbf{F}(M,N). 
\]
The induced homomorphism $i_*$ of zero homotopy groups is not surjective; even if the restriction of a fold map $f$ to the singular set $S(f)$ is an embedding, the map $f$ may not have a structure of a marked fold map (see Remark~\ref{r:9.5}). By Remark~\ref{r:6.1}, the induced homomorphism $i_*$ is not injective. Thus, $i_*$ is neither surjective, nor injective. 
\end{remark}

\section{The sheaf $\mathcal{F}$ of marked fold maps}\label{se:10}

Let $\mathcal{E}$ denote the category of smooth manifolds without boundary and smooth submersions. Then a set valued contravariant functor $F$ on $\mathcal{E}$ is said to be a {\it sheaf} on $\mathcal{E}$ if for any open covering $\{U_i\}$ of any manifold $X\in \mathcal{E}$, and any sections $s_i\in F(U_i)$ with $s_i=s_j$ over $U_i\cap U_j$ for all $i,j$, there is a unique section $s\in F(X)$ that restricts to $s_i$ over $U_i$ for each $i$. Given a closed subset $A$ of a manifold $X\in \mathcal{E}$ and a germ $s$ of sections of $F$ near $A$, the set $F(X, A; s)\subset F(X)$ is defined to be the subset of sections whose germ over $A$ coincides with $s$. 

\begin{remark} Our definition of a sheaf differs from the one in the papers \cite{MW} of Madsen and Weiss, and \cite{GMTW} by Galatius, Madsen, Tillmann and Weiss as we define a sheaf on the category $\mathcal{E}$ where morphisms are not smooth maps, but smooth submersions. Such a generalization is necessary as the sheaf of marked fold maps defined below is not a sheaf on the category of smooth manifolds without boundary and smooth maps.    
\end{remark}

These definitions also make sense in the case where the functor $F$ takes values in the category $\mathbf{Top}$ of topological spaces and continuous maps. 

\begin{definition} Let $F$ be a sheaf and $X$ a smooth manifold without boundary. Let $p: X\times \R\to \R$ be the projection onto the second factor. Two sections $s_0$ and $s_1$ of $F(X)$ are said to be \emph{concordant} if there is a section $s$ of $F(X\times \R)$ that agrees with the section $p^*s_0$ over $X\times (-\infty, -1)$ and with the section $p^*s_1$ over $X\times (1, \infty)$.  
\end{definition}

Let $\Delta_e$ denote the category with objects given by ordered sets $\underline{n}=\{0, 1, \dots, n\}$, one for each $n\ge 0$, and morphisms $\underline{m}\to \underline{n}$ given by order preserving maps. Every morphism in $\Delta_e$ can be written in terms of morphisms
\[
    \delta_i\colon \underline{n-1}\longrightarrow \underline{n}, \quad  \delta_i\colon \{0,1,\dots, n-1\}\mapsto \{0,1,\dots , i-1, i+1, \dots, n\}, 
\]
\[
    \sigma_i\colon \underline{n}\longrightarrow \underline{n-1}, \quad  \sigma_i\colon \{0,1,\dots ,n\}\mapsto \{0,1,\dots ,i-1, i, i, i+1,\dots, n-1\},
\]
called \emph{face} and \emph{degeneracy} maps respectively. This morphisms should be compared with their linear namesakes 
\[
\delta_i\colon \Delta^{n-1}\longrightarrow \Delta^n, \quad \delta_i(t_0, . . . , t_{n-1}) = (t_0, . . . , t_{i-1}, 0, t_i, . . . , t_{n-1}),
\]
\[
\sigma_i\colon \Delta^{n}\longrightarrow \Delta^{n-1}, \quad \sigma_i(t_0, . . . , t_{n}) = (t_0, . . . , t_i + t_{i+1}, . . . , t_{n})
\]
of standard geometric $n$-simplices
\[
\Delta^n = \{(t_0, \dots , t_n)\subset \R^{n+1} : 0 \le t_i \le 1, t_0+\cdots t_n=1\}.
\]
 A \emph{simplicial object} in a category $\mathcal{C}$ is a contravariant functor $F: \Delta_e\to \mathcal{C}$. A \emph{simplicial set} is a simplicial object in the category of sets, while a \emph{simplicial space} is a simplicial object in the category of topological spaces.   
 
Let $X$ be a simplicial space. Put $\partial_i=X(\delta_i)$ and $s_i=X(\sigma_i)$. Then the \emph{geometric realization} of a simplicial space $X$ is defined to be the topological space
\[
      |X|\colon=(\bigsqcup_{n\ge 0} X(\underline{n})\times \Delta^n)/\sim
\]
where the equivalence relation $\sim$ identifies $(x, \delta_i(u))$ with $(\partial_i(x), u)$ and $(x, \sigma_i(u))$ with $(s_i(x), u)$. 

For a manifold $U\in \mathcal{E}$, we say that a submanifold $W\subset U\times \R^{\infty}$ has a \emph{structure of an $\mathcal{F}$-submanifold} (or, simply, that $W$ is an \emph{$\mathcal{F}$-submanifold}) if the inclusion $W\subset U\times \R^{\infty}$ is the product of a marked fold map $W\to U$ and an arbitrary map $W\to \R^{\infty}$. Let $\mathcal{F}(U)$ denote the topological space of $\mathcal{F}$-submanifolds in $U\times \R^{\infty}$. Then the correspondence $U\mapsto \mathcal{F}(U)$ is a sheaf $\mathcal{F}$ called the {\it sheaf of marked fold maps}.

An extended $n$-simplex $\Delta^n_e$, with $n\ge 0$, is defined by
\[
    \Delta^n_e\colon=\{\ (x_0, ..., x_n)\in \R^{n+1}\ |\ x_0+ \cdots +x_n=1\}.
\]  
We note that the correspondence $\underline{n}\mapsto \mathcal{F}(\Delta_e^n)$ is not a simplicial object since for a morphism $\underline{m}\to \underline{n}$ the map $\mathcal{F}(\Delta_e^n)\to \mathcal{F}(\Delta_e^m)$ may not exist. However $\mathcal{F}$ do determine a simplicial set, which we describe next.

Let $\pi_W$ denote the projection $W\to \Delta^n_e$ of an $\mathcal{F}$-submanifold in $\Delta^n_e\times \R^{\infty}$ onto the first factor. In particular $\pi_W$ is a marked fold map. Let $\Sigma\subset W$ be the set of singular points of $\pi_W$. Then $\Sigma$ is a submanifold of $W$ and $\pi_W|\Sigma$ is an immersion. Let $S_k\subset \Delta_e^n$ denote the set of points $x$ such that the fiber of $\pi_W$ over $x$ is a singular manifold with precisely $k$ singular points. We will form a simplicial set by considering only those $\mathcal{F}$-submanifolds that satisfy properties $\mathbf{P1}$ and $\mathbf{P2}$. 

\begin{itemize}
\item[($\mathbf{P1}$)] For each morphism $\underline{m}\to\underline{n}$ given by a composition of face maps, the corresponding inclusion $\Delta^m_e\to \Delta^n_e$ is transversal to the immersed submanifolds $S_i$ for all $i$.  
\end{itemize}

The property ($\mathbf{P1}$) insures that for each morphism $\underline{m}\to \underline{n}$ the map $i: \Delta^m_e\to \Delta^n_e$ pulls back the map $\pi_W$ to a smooth map $i^*\pi_W: i^*W\to \Delta^m_e$. We note that the subspace in $\mathcal{F}(\Delta^n_e)$ of $\mathcal{F}$-submanifolds satisfying ($\mathbf{P1}$) is open dense. Here we consider only morphisms given by a composition of face maps since for degeneracy maps the corresponding statement is true ($\mathcal F$ is a sheaf). 

\begin{remark} The property ($\mathbf{P1}$) guarantees that $i^*\pi_W$ is a fold map. 
\end{remark}

\begin{remark}\label{r:9.2} We note that the subspace of $\mathcal{F}(\Delta_e^n)$ of $\mathcal{F}$-submanifolds satisfying the property ($\mathbf{P1}$) is not weakly homotopy equivalent to the space $\mathcal{F}(\Delta_e^n)$, since for example these two spaces have different number of path components.  
\end{remark}

\begin{itemize}
	\item[($\mathbf{P2}$)] The map $i^*\pi_W$ is a marked fold map with the structure inherited from the marked fold map $\pi_W$.  
\end{itemize}

\begin{remark} Let us describe the property ($\mathbf{P2}$) in terms of markings. Let $W$ be an $\mathcal{F}$-submanifold in $\mathcal{F}(\Delta_e^n)$ satisfying the property ($\mathbf{P1}$) and let $S$ be a path component of $S_k$. Let 
\[
   i_{\Delta}\colon D^k\times S \longrightarrow \Delta_e^n, 
\]
\[
   i_W\colon V\times S\longrightarrow W
\]
be local marking maps, where $V$ is the source manifold of an elementary $k$-bordism. The map $i^*\pi_W$ inherits a local marking if there are commutative diagrams:
\[
\begin{CD}
   D^k\times (S\cap \Delta^m_e) @>>> \Delta^m_e \\
   @VVV @VVV \\
   D^k\times S @>>> \Delta^n_e
\end{CD}
\]
and
\[
\begin{CD}
   V\times (S\cap \Delta^m_e) @>>> W \\
   @VVV @VVV \\
   V\times S @>>> W, 
\end{CD}
\]
where the vertical maps are inclusions. In particular, the property ($\mathbf{P2}$) requires the existence of the two commutative diagrams above. 

Similarly, let $\Sigma$ be a path component of fold singular points in $W$ and let 
\[
    \alpha\colon U\times \Sigma \longrightarrow W \quad\mathrm{and}\quad \beta\colon [-1,1]\times \Sigma \longrightarrow \Delta^n_e
\]
be global marking maps, where $U$ is a standard handle. Then, the map $i^*\pi_W$ inherits the global marking if there are commutative diagrams
\[
\begin{CD}
   U\times (\Sigma\cap i^*W) @>>> i^*W \\
   @VVV @VVV \\
   U\times \Sigma @>>> W
\end{CD}
\]
and
\[
\begin{CD}
   [-1,1]\times (\Sigma\cap i^*W) @>>> \Delta_e^m \\
   @VVV @VVV \\
   [-1,1]\times \Sigma @>>> \Delta_e^n, 
\end{CD}
\]
where the vertical maps are inclusions. Again, the property ($\mathbf{P2}$) requires the existence of the two commutative diagrams above. 

It is easily verified that if the local and global markings for $\pi_W$ are compatible, then the same is true for inherited markings for $i^*\pi_W$. 
\end{remark}

\begin{remark}
The space of $\mathcal{F}$-submanifolds that satisfy ($\mathbf{P1}$)
is weakly homotopy equivalent to the space of $\mathcal{F}$-submanifolds that satisfy ($\mathbf{P1}$) and ($\mathbf{P2}$). 
\end{remark}

Let $\mathcal{F}^{\pitchfork}(\Delta^n_e)\subset \mathcal{F}(\Delta^n_e)$ denotes the set  of those $\mathcal{F}$-submanifolds $W$ in $\Delta^n_e\times \R^{\infty}$ that satisfy the properties $(\mathbf{P1})$ and $(\mathbf{P2})$. 
Then the correspondence $\underline{n}\mapsto \mathcal{F}^{\pitchfork}(\Delta^n_e)$ is a simplicial object. For a morphism $i: \underline{m}\to \underline{n}$, the map 
\[
   \mathcal{F}^{\pitchfork}(i)\colon \mathcal{F}^{\pitchfork}(\underline{n})\longrightarrow \mathcal{F}^{\pitchfork}(\underline{m})
\]
takes an $\mathcal{F}$-submanifold $W$ in $\Delta_e^n\times \R^{\infty}$ to the $\mathcal{F}$-submanifold $\tilde{i}^*W$ in $\Delta^m_e\times \R^{\infty}$, where $\tilde{i}$ is the linear map $\Delta^m_e\to \Delta^n_e$ corresponding to the morphism $i$. 

The realization of the simplicial object $\underline{n}\mapsto \mathcal{F}^{\pitchfork}(\Delta^n_e)$ is denoted by $|\mathcal{F}|$.

\section{Classifying space $|\mathcal{F}|$ of the sheaf $\mathcal{F}$}\label{se:11}

Proposition~\ref{p:4} will show that the space $|\mathcal{F}|$ is a classifying space of the sheaf $\mathcal{F}$. In view of Remark~\ref{r:9.2} its statement is slightly surprising. However the proof of Proposition~\ref{p:4} is similar to that of \cite[Proposition 2.4.3]{MW}.

\begin{proposition}\label{p:4} Let $z$ be a value in $\mathcal{F}(pt)$. Then for any manifold $X$ in $\mathcal{E}$ and a closed subset $A\subset X$, there is a natural bijection of the set of homotopy classes of maps $(X, A)\to (|\mathcal{F}|, z)$ and $\mathcal{F}[X, A; z]$.  
\end{proposition}
\begin{proof} As has been mentioned the proof of Proposition~\ref{p:4} follows the proof of \cite[Proposition 2.4.3]{MW}. In contrast to the functor considered in \cite{MW}, the functor $\mathcal{F}$ is defined not on the category of smooth manifolds and smooth maps, but on the category of smooth manifolds and smooth embeddings. In particular, a smooth map $X\to X$ does not induce a morphism $\mathcal{F}(X)\to \mathcal{F}(X)$ used in the proof of \cite[Proposition 2.4.3]{MW}. Because of that to prove Proposition~\ref{p:4} a slight modification of the argument in the proof of \cite[Proposition 2.4.3]{MW} is necessary. We will explain the necessary modification in the case where $A$ is empty; the  general case is similar. 

Let $X$ be a smooth manifold in $\mathcal{E}$ and $u$ an element in $\mathcal{F}(X)$ representing a class $[u]$ in $\mathcal{F}[X]$, i.e., $u$ is given by a submanifold $W\subset X\times \R^{\infty}$ such that its projection $\pi_X: W\to X$ onto the first factor is a proper marked fold map. There is a smooth triangulation of $X$ with an ordered set of vertices. Following the proof of \cite[Proposition 2.4.3]{MW}, for each characteristic embedding $c: \Delta^n\to X$ of the triangulation we consider an embedding $c_e: \Delta^n_e\to X$ of an extended simplex. If necessary we may modify the representative $u$ in the class $[u]$ and the triangulation so that 

\begin{itemize}
\item[($\mathbf{C1}$)] for each extended characteristic embedding $c_e$, the pullback $c_e^*u$ is an element in $\mathcal{F}^{\pitchfork}(\Delta^n_e)$. 
\end{itemize}

Indeed, first, we may modify the triangulation so that the characteristic embeddings $c_e$ are transversal to the submanifolds $S_k$ of $X$ of points parametrizing fibers of $\pi_X$ with $k$ singular points. Next, we may modify the markings of the fold map $\pi_X$ so that the image of each marking map is in a small tubular neighborhood of the set of singular points or singular values. Finally we may modify the triangulation so that the condition ($\mathbf{C1}$) is satisfied. 

If the condition ($\mathbf{C1}$) is satisfied, then each map $c_e$ determines a smooth map from $\mathop{\mathrm{Im}}(c_e|\Delta^n)\subset X$ to a simplex in $|\mathcal{F}|$. Thus there is a well defined map $\psi: X\to |\mathcal{F}|$. 

Given a different smooth triangulation of $X$ with an ordered set of vertices, there is a triangulation of $X\times [0,1]$ with an ordered set of vertices such that the restrictions of the triangulation to $X\times \{0\}$ and $X\times \{1\}$ coincide with the former and latter triangulations of $X$ respectively. Furthermore, we may assume that the triangulation of $X\times [0,1]$ satisfies the condition ($\mathbf{C1}$). Then the obtained map $X\times [0,1]\to |\mathcal{F}|$ is a homotopy of two maps $X\to |\mathcal{F}|$, which shows that the homotopy class of $\psi$ does not depend on the triangulation of $X$.   

Let now $g: X\to |\mathcal{F}|$ be a map representing a class $[g]\in [X, |\mathcal{F}|]$. By the simplicial approximation principle (see the proof of \cite[Proposition 2.4.3]{MW}), there is a simplicial set $X^s$ and a simplicial map of $X^s$ into the simplicial set $\underline{n}\mapsto \mathcal{F}^{\pitchfork}(\Delta^n_e)$ whose realization is $g$. There is a smooth homotopy $h_t: X\to X$ of $h_0=\id_X$, with $t\in [0,1]$ such that 
\begin{itemize}
\item $h_t$ maps each simplex of $X$ onto itself for each $t$, and
\item each simplex has a neighborhood that maps by $h_1$ onto the simplex.
\end{itemize}
We may also assume that
\begin{itemize}
	\item for each simplex $\Delta$ in $X$, the restriction of $h_t$ to $h_t^{-1}(\mathring{\Delta})$ is a submersion to the interior $\mathring{\Delta}$ of $\Delta$, i.e., at each point in $h^{-1}(\mathring{\Delta})$, its rank is $\dim \Delta$.  
\end{itemize}
Our additional assumption on the homotopy $h_t$ guarantees that the element constructed in the proof of \cite[Proposition 2.4.3]{MW} is in $\mathcal{F}(X)$. 
\end{proof}

\section{Equivalence $|\mathcal{F}|=\mathbf{T}$}\label{se:12}

Let $f_i\colon M_i\to N_i$ be two smooth maps of oriented manifolds, $i=1,2$. Without loss of generality we may assume that 
\[ 
   r=\dim N_2 - \dim N_1
\]
is non-negative. We say that the fiber of $f_1$ over a point $y_1\in N_1$ is \emph{equivalent} to the fiber of $f_2$ over a point $y_2\in N_2$ if there are open neighborhoods $U_1$ of $y_1$ in $N_1$ and $U_2$ of $y_2$ in $N_2$ and there are orientation preserving diffeomorphisms 
\[
    \alpha\colon U_1\times \R^r\longrightarrow U_2 
\]
and
\[
   \beta\colon f^{-1}_1(U_1)\times \R^r\longrightarrow f_2^{-1}(U_2)
\]
for which the diagram
\[
\begin{CD}
f_1^{-1}(U_1)\times \R^r @>\beta>> f_2^{-1}(U_2) \\
@Vf_1 VV @Vf_2 VV \\
U_1\times \R^r @>\alpha>> U_2
\end{CD}
\]
commutes. Similarly we can define equivalence classes of fibers of marked fold maps. By Remark~\ref{rem:9.5} an equivalence class of fibers of marked fold maps determines an equivalence class in $\I$ of an oriented compact manifold with parametrized boundary. 

Each point $x$ in $|\mathcal{F}|$ belongs to the interior of a unique non-degenerate simplex, which corresponds to an $\mathcal{F}$-submanifold $W$ in $\Delta^n_e\times \R^{\infty}$. In particular each point $x$ determines the equivalence class of the fiber over $x$ under the projection $\pi: W\to \Delta^n_e$, which is a marked fold map. Hence, there is a natural stratification 
\[
     |\mathcal{F}|=\bigsqcup_{\alpha\in \I} |\mathcal{F}|_{\alpha}.
\]

For each $u\in \mathcal{F}^{\pitchfork}(\Delta^n_e)$, let $W_u$ be the $\mathcal{F}$-submanifold in $\Delta^n\times \R^{\infty}$, $n=0,1,...$. Then, for each $i$, the map 
\[
   \mathcal{F}^{\pitchfork}(\Delta^{n-1}_e)\times \Delta^{n-1}\longrightarrow \mathcal{F}^{\pitchfork}(\Delta^n_e)\times \Delta^n
\]
that takes $(\partial_i(x), u)$ to $(x, \delta_i(u))$ is covered by the map

\begin{equation}\label{e:11.1}
   \bigsqcup_{u\in \mathcal{F}^{\pitchfork}(\Delta^{n-1}_e)} W_u \longrightarrow \bigsqcup_{v\in \mathcal{F}^{\pitchfork}(\Delta^n_e)} W_v
\end{equation}
and the map 
\[
   \mathcal{F}^{\pitchfork}(\Delta^{n}_e)\times \Delta^{n}\longrightarrow \mathcal{F}^{\pitchfork}(\Delta^{n-1}_e)\times \Delta^{n-1}
\]
that takes $(s_i(x), u)$ to $(x, \sigma_i(u))$ is covered by the map

\begin{equation}\label{e:11.2}
   \bigsqcup_{u\in \mathcal{F}^{\pitchfork}(\Delta^{n}_e)} W_u \longrightarrow \bigsqcup_{v\in \mathcal{F}^{\pitchfork}(\Delta^{n-1}_e)} W_v.
\end{equation}

The disjoint union $\sqcup W_u$ of all $\mathcal{F}$-submanifolds $W_u$ subject to identifications (\ref{e:11.1})
and (\ref{e:11.2}) is denoted by $\mathbf{W}$. There is a canonical map 
\[
    \Pi: \mathbf{W} \longrightarrow |\mathcal{F}|
\]
that restricts to the projection of $\mathcal{F}$-submanifolds $W_u$ to $\Delta^n$, each of the projections is a marked fold map.

\begin{proposition}\label{pr:9.5} For each $\alpha\in \I$, the space $|\mathcal{F}|_{\alpha}$ is weakly homotopy equivalent to the space $\BDiff M_{\alpha}$. 
\end{proposition}
\begin{proof} Suppose that $\alpha$ is in $\I_k$. Let $y$ be a point in $|\mathcal{F}|_{\alpha}$. It belongs to the interior of a unique non-degenerate simplex in  $|\mathcal{F}|$ and therefore corresponds to a marked fold map $\pi: W\to \Delta^n_e$. The set $|\mathcal{F}|_{\alpha}$ has a trivial normal bundle $U$ in $|\mathcal{F}|$ of dimension $k$, the fiber $U_y$ of $U$ over $y$ is trivialized by the local marking of the marked fold maps $W\to \Delta^n_e$. The restriction of $\Pi$ to $\Pi^{-1}(U)$ will be denoted by $\Pi_{\alpha}: \mathbf{W}_{\alpha}\to U$.  

In fact we may choose $U$ so that for each $y$ the set $W_y:=\pi^{-1}(U_y)$ is a submanifold of $W$ and $\pi|W_y: W_y\to U_y$ is a $k$-bordism with a kernel. Then $\Pi_{\alpha}$ is a family of $k$-bordisms $\pi|W_y$ with kernels parametrized by $|\mathcal{F}|$. Note that the complement in $W_y$ to the image of the kernel of the $k$-bordism $\pi|W_y$ is the total space of a smooth bundle over $U_y$ with fiber $M_{\alpha}$. Thus, by removing images of kernels of $k$-bordisms, we obtain a bundle over $U$ with fiber $M_{\alpha}$. Its structure group is $\Diff M_{\alpha}$ since the boundary of each fiber comes with a parametrization. Hence this bundle is classified by a map

%\begin{remark}\label{rem:12.1} For each family $f: W\times P\to X\times P$ of $k$-bordisms with kernels parametrized %by a manifold $P$, the construction in the proof of Proposition \ref{p:4} gives a map $X\times P\to |\mathcal{F}|$. In %fact, the proof of Proposition \ref{p:4} shows that not only the homotopy class of the map $X\times P\to %|\mathcal{F}|$ is well defined, but the same is true for the homotopy class of $X\times P\to U$.  
%\end{remark}
\[ 
     \psi\colon U\longrightarrow \BDiff M_{\alpha}.  
\] 
Since we only need to show that $|\mathcal{F}|_{\alpha}$ is weakly homotopy equivalent to the space $\BDiff M_{\alpha}$, we may replace $\BDiff M_{\alpha}$ by a CW complex weakly homotopy equivalent to $\BDiff M_{\alpha}$. We will continue to denote the new CW complex by $\BDiff M_{\alpha}$.  Then we may use the universal bundle $\mathbf{E}_0\to \BDiff M_{\alpha}$ in order to form a family of $k$-bordisms with kernels over $\BDiff M_{\alpha}\times D^k$. For each point $x\in \BDiff M_{\alpha}$, the $k$-bordism with kernel over $x\times D^k$ is classified by a map to $|\mathcal{F}|$ with image in $U$. The union of these maps determine a map  
\[
    \varphi\colon \BDiff M_{\alpha}=\BDiff M_{\alpha}\times\{0\}\subset \BDiff M_{\alpha}\times D^k\to U.
\] 
Clearly $\psi$ and $\varphi$ are inverses of each other. 
\end{proof}

Thanks to Proposition~\ref{pr:9.5}, the space $|\mathcal{F}|$ (or, more precisely, a space weakly homotopy equivalent to $|\mathcal{F}|$) has a fairly simple description, which we give in this subsection. 

Let $\I_k\subset \I$ denote the subset indexing diffeomorphism types of oriented compact manifolds with precisely $k$ boundary components of the form $S^{i-1}\times S^{d-i}$. Let $\D_0$ be a CW complex equipped with a weak homotopy equivalence
\[
    d_0: \D_0\colon = \bigsqcup_{\alpha\in \mathfrak{J}_0} \D(M_{\alpha}) \longrightarrow \bigsqcup_{\alpha\in \mathfrak{J}_0} \BDiff M_{\alpha}.
\]
Then $\D_0$ is a space classifying oriented manifold bundles with fiber of dimension $d$; the universal oriented bundle $\E_0\to \D_0$ is given by the  pullback via $d_0$ of the universal oriented manifold bundle over $\sqcup \BDiff M_{\alpha}$. 

\begin{remark} We note that the universal oriented bundle $\E_0\to \D_0$ has no smooth structure. In particular, for a continuous map $\tau: N\to \D_0$, the total space $\tau^*\E_0$ of the pullback has no structure of a smooth manifold. However it is possible to introduce a notion of a smooth map of a manifold into $\D_0$ so that for a smooth map $\tau$ to $\D_0$, the pullback bundle is smooth.
\end{remark}

Suppose that there is a universal space $\D_{k-1}$ for families of marked fold maps with fibers with at most $k-1$ singular points such that $\D_{k-1}$ is a CW complex. In other words $\D_{k-1}$ is a CW complex such that families of marked fold maps into a space $N$ over $X$ with fibers with at most $k-1$ singular points are in bijective correspondence with homotopy classes of maps $[N, \D_{k-1}]$. 

Let $\D'_k$ be a CW complex equipped with a weak homotopy equivalence
\[
    \D'_{k}\colon= \bigsqcup_{\alpha\in\I_k} \D(M_{\alpha})\longrightarrow \bigsqcup_{\alpha\in\I_k} \BDiff M_{\alpha},
\]
where $\BDiff M_{\alpha}$ stands for the group of diffeomorphisms of $M_{\alpha}$ trivial on $\partial M_{\alpha}$.

Then there is a family $\E'\to \D'_k\times D^k$ constructed as in Definition~\ref{d:8} such that its restriction to $\D'_k\times \partial D^k$ is a special family of order $k$ (see Definition~\ref{d:8}). It is classified by a map into $\D_{k-1}$. Then the homotopy pushout $\D_k$ of
\[
      \D'_k\times D^k\longleftarrow \D'_k\times \partial D^k\longrightarrow \D_{k-1} 
\]
is a classifying space for marked fold maps with fibers with at most $k$ singular points. 

\begin{remark} The proof of Proposition~\ref{pr:9.5} implies that there is a map $\mathbf{E}_k\to \D_k$, which, in some certain sense, is a universal marked fold map with at most $k$-fold fold points. 
\end{remark}

\begin{remark} Our argument in the construction of $\D_{k}$ is justified by Proposition~\ref{pr:9.5}; the construction of $\D$ is only a description of the space $|\mathcal{F}|$ according to Proposition~\ref{pr:9.5}. 
\end{remark}

\begin{remark}
It is tempting to use an induction with inductive assumption that there exists a universal map $u_k: \mathbf{E}_{k-1}\to \D_{k-1}$ for marked fold maps with fibers with at most $k-1$ singular points such that $\D_{k-1}$ is a CW complex and for each marked fold map $f: M\to N$ with fibers with at most $k-1$ singular points, there is a unique up to homotopy map $g: N\to \D_{k-1}$ such that $f=g^*u_k$. However, in order to establish the induction step, one needs to
\begin{itemize}
\item enlarge the class of fold maps so that the new class includes special families of order $k$, and  
\item devise a consistent method of how to introduce smooth structures on pullbacks $g^*\mathbf{E}_{k-1}$ for maps $g: N\to \D_{k-1}$. 
\end{itemize}

The latter is not immediate since the space $\mathbf{E}_{k-1}$ has no smooth structure. One should address this question with care, since, for example, there exist two non fiberwise diffeomorphic smooth fiber bundles that are isomorphic as topological fiber bundles (see the Hatcher example in \cite{Ig}).   

\end{remark}

The construction proceeds by induction. 

It immediately follows that the space $\D$ classifies marked fold maps in the sense that for each oriented closed manifold $N$, the concordance classes of marked fold maps into $N$ are in bijective correspondence with elements in $[N, \D]$ (see Definition~\ref{def:20.1}).

\section{Equivalence $\D\approx \mathfrak{H}'(d)$} 
\subsection{The space $\D$ as a direct limit.} In this subsection we represent the classifying space $\D$ of marked fold maps as a direct limit of spaces $\c(M_{\alpha})$. 

Let $\I_k(\alpha)$ denote the set of indexes $\beta\in \I_k$ for which the manifold with paramet\-rized boundary $M_{\alpha}$ can be obtained by attaching handles to $M_{\beta}$ along some of the boundary components of $\partial M_{\beta}$. Clearly $\I_k(\alpha)$ is empty for $k<l$ if $\alpha\in \I_l$. 
Let
\[
     f_{\beta}: W\longrightarrow D^k  
\]  
be the elementary $k$-bordism associated with $M_{\beta}$ by means of Lemma~\ref{l:4.1}. Let $\mathfrak{I}(\beta, \alpha)$ be the set indexing the path components of the subset of $D^k$ of points paramet\-rizing fibers of $f_{\beta}$ equivalent to $M_{\alpha}$. In particular, the set $\I(\beta, \alpha)$ is empty unless $\beta\in \I_k(\alpha)$. For $\kappa\in \mathfrak{I}(\beta, \alpha)$ let $D^k_{\kappa}$ denote the closure in $D^k$ of the set indexed by $\kappa$. 

\begin{example} Let $M_{\alpha_0}$ be a compact oriented surface of genus $g$ with one boundary component of the form $S^0\times S^1$. Then the set $\I_3(\alpha_0)$ contains precisely three indices $\alpha_0$, $\alpha_1$ and $\alpha_2$ corresponding to oriented compact surfaces $M_{\alpha_0}$, $M_{\alpha_1}$ and $M_{\alpha_2}$ of genus $g$ with respectively one, two and three boundary components of the form $S^{0}\times S^{1}$. For $i=0,1,2$, the set $\sqcup D_{\kappa}^i$ where $\kappa$ ranges over $\I(\alpha_i, \alpha_0)$ respectively is the point, the \underline{disjoint union} $[-1,0]\sqcup [0,1]$, and the \underline{disjoint union} of the four closed quarters of a disc of dimension $2$. 
\end{example}

To begin with we construct spaces $\c(M_{\alpha})$, one for each $\alpha\in \I$, by induction similar to that in the construction of $\D$. In the construction of $\c(M_{\alpha})$ the role of $\D'_k$ is played by a CW complex $\c'_k(M_{\alpha})$ equipped with a weak homotopy equivalence 
\[
       \c'_k(M_{\alpha})\colon = \bigsqcup_{\substack{\beta\in \I_k(\alpha)\\ \kappa\in \I(\beta, \alpha)}} \c'(M_{\beta}, M_{\alpha}; \kappa)\longrightarrow \bigsqcup_{\beta\in \I_k(\alpha)} \BDiff M_{\beta}\times \I(\beta, \alpha), 
\]
where $\c'(M_{\beta}, M_{\alpha}; \kappa)=\BDiff M_{\beta}\times \{\kappa\}$ is a CW complex weakly homotopy equivalent to $\BDiff M_{\beta}$. 
%\[
%       \c'_k(M_{\alpha})\colon = \bigsqcup_{\beta\in \I_k(\alpha)} %\c'(M_{\beta}, M_{\alpha}), 
%\]
%each component $\c'(M_{\beta}, M_{\alpha})$ of which is the total space of a %covering  
%\begin{equation}\label{eq:7.2}
%      \pi_c:  \c'(M_{\beta}, M_{\alpha}) \longrightarrow  \D(M_{\beta})
%\end{equation}
%with fiber $\I(\beta, \alpha)$. To describe the monodromy of the covering %(\ref{eq:7.2}), let us observe that each element $l$ in %$\pi_1(\D(M_{\beta}))$ corresponds to an equivalence of $(M_{\beta}, %w_{\beta})$, and therefore, to an equivalence of $(W, f_{\beta})$. In %particular, each element $l$ gives rise to a permutation of factors of $D^k$, %which in its turn determines a permutation of elements in the fiber %$\I(\beta,\alpha)$ of $\pi_c$.

The CW complex $\c_0(M_{\alpha})$ is defined to be $\BDiff M_{\alpha}$. For $k\ge 0$, the CW complex $\c_k(M_{\alpha})$ is defined to be the pushout
\[
     \bigsqcup_{\substack{\beta\in \I_k(\alpha)\\ \kappa\in \I(\beta, \alpha)}} \c'(M_{\beta}, M_{\alpha}; \kappa) \times D^k_{\kappa} \longleftarrow      \bigsqcup_{\substack{\beta\in \I_k(\alpha)\\ \kappa\in \I(\beta, \alpha)}} \c'(M_{\beta}, M_{\alpha}; \kappa)\times  (D^k_{\kappa}\cap \partial D^k)\longrightarrow \c_{k-1}(M_{\alpha}). 
\]
In fact we may arrange the maps so that there is an obvious canonical map 
\[
\eta_{k, \alpha}\colon \c_k(M_{\alpha})\to \D_k.
\]  

We note that if $\alpha\in \I_l$, then $\c'_k(M_{\alpha})$ is empty for $k<l$. For $k=l$, the map $\eta_{k, \alpha}$ is a homeomorphism onto image. If $k>l$, then each point $x\in \D_k$ has at most finitely many inverse images in $\c_k(M_{\alpha})$. 

Finally the space $\c(M_{\alpha})$ is defined to be the colimit of spaces $\c_k(M_{\alpha})$ with respect to inclusions $\c_k(M_{\alpha})\subset \c_{k+1}(M_{\alpha})$. Under the above assumptions on choices of maps, there is a canonical map 
\[
    \eta_{\alpha}\colon \c(M_{\alpha})\longrightarrow \D
\]
for each $\alpha\in \I$. These maps commute with obvious inclusions 
\begin{equation}\label{eq:12.1}
   \c(M_{\beta})\longrightarrow \c(M_{\alpha}). 
\end{equation}

\begin{remark} The space $\D$ is identical to the direct limit of spaces $\c(M_{\alpha})$ taken with respect to obvious inclusions (\ref{eq:12.1}).
\end{remark}

\subsection{Homotopy colimit}

Let us now define a new functor $\mathcal{D}'$ from $\I$ to the category of topological spaces. By definition, $\mathcal{D}'(\alpha)=\c(M_{\alpha})$ for each $\alpha\in \I$, where $M_{\alpha}$ is a manifold with parametrized boundary representing the equivalence class of $\alpha$. To each arrow $\beta\to \alpha$, the functor $\mathcal{D}'$ assigns an obvious inclusion 
\[
   \c(M_{\beta})\longrightarrow \c(M_{\alpha}). 
\]  
Put
\[
      \mathfrak{H}'(d)\colon = \mathop{\mathrm{hocolim}} \c(M_{\alpha}),
\]
where the homotopy colimit is taken with respect to the diagram $\mathcal{D}'$.

For the next theorem we need a complex $\Delta_k$. It is associated with the 
configuration of coordinate planes in $\R^k$. To begin with we construct an oriented graph $\Gamma$. Its vertices are in bijective correspondence with path components of half planes given by $k$ equations and inequalities
\[
    x_i\ge 0, \quad x_i\le 0 \quad \mathrm{or} \quad x_i=0,
\]
one for each coordinate $x_i$ in $\R^k$ so that $\Delta_k$ has precisely $3^k$ vertices. 
If a half plane $L_0$ of dimension $s$ is in the boundary of a half plane $L_1$ of dimension $>s$, then there is an oriented edge from the vertex of $L_0$ to that of $L_1$.  
The complex $\Delta_k$ is constructed so that its simplices of dimension $l$ correspond to paths in $\Gamma$ of length $l+1$. These simplices are glued along the boundaries in an obvious way. We note that the space $\Delta_k$
is contractible as there is a deformation of $\Delta_k$ that linearly deforms $(k+1)$-simplices of $\Delta_k$ to the unique common point.

Let us recall that a continuous map $p: E\to B$ to a path connected space $B$ is said to be a \emph{quasifibration} if the induced map
\[
     p_*: \pi_i(E, p^{-1}(y_0), x_0) \longrightarrow \pi_i(B, y_0)
\]
is an isomorphism for each $y_0\in B$, $x_0\in \pi^{-1}(y_0)$ and $i\ge 0$. In fact a map $p: E\to B$ is a quasifibration if and only if the base space $B$ is a union of open subsets $V_1$ and $V_2$ such that the restrictions of $p$ to $V_1$, $V_2$ and $V_1\cap V_2$ are quasifibrations (e.g., see \cite{Hat1}). Another criterion asserts that $p: E\to B$ is a quasifibration if there are deformation retracts $E'$ of $E$ and $B'$ of $B$ under deformations $F_t$ and $f_t$ respectively with $t\in [0,1]$ such that $F_t$ covers $f_t$, $E'\to B'$ is a quasifibration, and $F_1: p^{-1}(x)\to p^{-1}(f_1(x))$ is a weak homotopy equivalence for each $x\in B$. 

\begin{theorem}\label{th:4.1}  The projection $p: \mathfrak{H}'(d)\to \D$ is a quasifibration. 
\end{theorem}

\begin{proof} To begin with let us describe the fiber of $p$ over a point $x$ in $\D_k\setminus \D_{k-1}$ for $k\ge 0$. Here we assume that $\D_{-1}$ is empty. 
%Suppose that $x$ corresponds to the class of an oriented compact manifold $M_{\alpha}$ with $\alpha\in \I_k$. 
Then the fiber $p^{-1}(x)$ is a CW complex of dimension $k$ isomorphic to $\Delta_k$. In fact, there is a canonical map from the fiber $p^{-1}(x)$ to $\Delta_k$. Hence, the restriction of $p$ to $p^{-1}(\D_k\setminus \D_{k-1})$ is a trivial fiber bundle with fiber $\Delta_k$. We note that $\D_{k-1}$ is open in $\D_k$ and that $\D_{k}\setminus \D_{k-1}$ has an open neighborhood in $\D_k$ which is a trivial disc bundle $(\D_k\setminus \D_{k-1})\times \mathop{\mathrm{Int}}(D^k)$.

Furthermore, the radial retraction  
\[
 U_k\colon=(\D_k\setminus \D_{k-1})\times \mathop{\mathrm{Int}}(D^k)\longrightarrow \D_k\setminus \D_{k-1}
\]
is covered by a retraction of $p^{-1}(U)$ to $p^{-1}(\D_k\setminus \D_{k-1})$ that takes a fiber $\Delta_l$ over a point in $\D_l\setminus \D_{l-1}$ to a fiber $\Delta_k$ in an obvious way. Consequently, the restriction of $p$ to an open neighborhood of $\D_k\setminus \D_{k-1}$ is a quasifibration. It is also clear that the restrictions of $p$ to $U_k\cap U_l$ are quasifibrations. Hence $p$ is a quasifibration. 
\end{proof}

Theorem~\ref{th:4.1} implies that the homotopy fiber of $p$ is homotopy trivial since the inclusion of any fiber to the homotopy fiber of a quasifibration is a weak homotopy equivalence. In particular, the map $p$ induces an isomorphism of homotopy groups. 

\begin{corollary} The projection $p: \mathfrak{H}'(d)\to \D$ is a weak homotopy equivalence. 
\end{corollary}

\section{Equivalence $\mathfrak{H}'(d)\approx \mathfrak{H}(d)$}\label{se:14}

We have defined two diagrams $\mathcal{D}$ and $\mathcal{D}'$ over the same category $\I$. In this section we show that their homotopy colimits are weakly homotopy equivalent. In fact our argument shows that for any choice of models of spaces $\BDiff M_{\alpha}$ and maps $\BDiff M_{\beta}\to \BDiff M_{\alpha}$, the corresponding homotopy colimit is weakly homotopy equivalent to $\mathfrak{H}'(d)$.

To begin with, let us observe that for each $\alpha\in \I$ there is a weak homotopy equivalence
\[
   i_{\alpha}\colon  \mathcal{D}'(\alpha)=\mathbf{C}(M_{\alpha})\approx \BDiff M_{\alpha}= \mathcal{D}(\alpha).
\]
Its homotopy class is canonical, but the map $i_{\alpha}$ itself is not canonical. Similarly, for each arrow $\beta\mapsto \alpha$ in $\I$ the diagram 
\[
\begin{CD}
   \mathcal{D}'(\beta) @>\mathcal{D}'(\beta\mapsto \alpha)>> \mathcal{D}'(\alpha) \\
   @V i_{\beta} VV @Vi_{\alpha} VV \\
   \mathcal{D}(\beta) @>\mathcal{D}(\beta\mapsto\alpha)>> \mathcal{D}(\alpha),
\end{CD}
\] 
is homotopy commutative, but not necessarily commutative. Since the diagrams over arrows in $\I$ are not necessarily commutative, we can not immediately conclude that the homotopy colimit $\mathfrak{H}'(d)$ of $\mathcal{D}'$ is weakly homotopy equivalent to the homotopy colimit $\mathfrak{H}(d)$ of $\mathcal{D}$. I learned the following example from David Ayala~\cite{Ay}. 

\begin{example}\label{ex:3.2} Let us consider a diagram $\mathcal{D}_1:$  
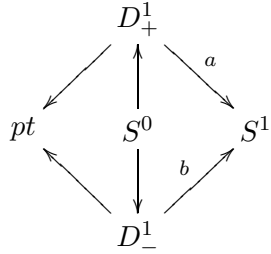
\[
\xymatrix
{
& & D^1_+ \ar[dl] \ar[dr]^{a}&  \\
& pt & S^0 \ar[u] \ar[d] & S^1 \\
& & D^1_- \ar[ul] \ar[ur]^b &
}
\]
where each of the two vertical maps is an inclusion of $S^0$ into the boundary of a copy of the segment $D^1$; and the maps $a, b$ are inclusions of the upper and lower semicircles into a circle. The homotopy colimit of the diagram $\mathcal{D}_1$ is contractible. 

Let now $\mathcal{D}_2$ be the diagram that differs from the diagram $\mathcal{D}_2$ only in that the maps $a$ and $b$ are constant maps onto a point in $S^1$. Then the homotopy colimit of $\mathcal{D}_2$ is homotopy equivalent to $S^1\vee S^2$.
Thus 
\[
   \mathop{\mathrm{hocolim}} \mathcal{D}_1\ncong \mathop{\mathrm{hocolim}}\mathcal{D}_2,
\]
though each map of the diagram $\mathcal{D}_2$ is homotopic to the corresponding map of the diagram $\mathcal{D}_1$.
\end{example}

Nevertheless we will show that the following theorem holds true. 

\begin{theorem}\label{th:14.1} The spaces $\mathfrak{H}(d)$ and $\mathfrak{H}'(d)$ are weakly homotopy equivalent.  
\end{theorem}

In order to establish Theorem~\ref{th:14.1} we will need two series of auxiliary homotopy colimits taken with respect to series of functors $\mathcal{D}_k$ and $\mathcal{D}'_k$. 

For each $k\ge 0$, the diagram $\mathcal{D}_k$ is defined to be the restriction of the functor $\mathcal{D}$ to the full subcategory $\I(k)\subset \I$ of indices $\alpha$ of diffeomorphism types of manifolds with at most $k$ boundary components. The homotopy colimit taken with respect to $\mathcal{D}_k$ is denoted by $\mathfrak{H}_k(d)$. 

The diagram $\mathcal{D}'_k$, with $k\ge 0$, is defined similarly to the diagram $\mathcal{D}'$. Namely, we put $\mathcal{D}'_k(\alpha)=\c_k(M_{\alpha})$ for each $\alpha\in \I$, where $M_{\alpha}$ is a manifold with parametrized boundary representing the equivalence class of $\alpha\in \I_k$. To each arrow $\beta\to \alpha$ in $\I_k$, the functor $\mathcal{D}'_k$ assigns an obvious inclusion 
\[
   \c_k(M_{\beta})\longrightarrow \c_k(M_{\alpha}). 
\]  
Put
\[
      \mathfrak{H}'_k(d)\colon = \mathop{\mathrm{hocolim}} \c_k(M_{\alpha}),
\]
where the homotopy colimit is taken with respect to the diagram $\mathcal{D}'_k$. We also recall that there are canonical maps 
\[
    \eta_{\alpha, k}\colon \mathcal{D}'_k(\alpha)=\c_k(M_{\alpha})\longrightarrow \mathbf{T}=|\mathcal{F}|. 
\]

\begin{remark} The functor $\mathcal{D}'_k$ is not a restriction of the functor $\mathcal{D}'$. Lemma~\ref{l:14.2} below is not true if we replace $\mathcal{D}'_k$ by the functor given by the restriction of $\mathcal{D}'$ to $\I_k$.  
\end{remark}

\begin{lemma}\label{l:14.2}
Let $x$ be a point in $\mathcal{D}'_k(\alpha)$. Then there is a unique index $\sigma$ and a point $s\in \mathcal{D}'_k(\sigma)$ such that for each point 
$y\in \mathcal{D}'_k(\gamma)$ that maps by $\mathcal{D}'_k(\gamma\mapsto \alpha)$ onto $x$, there is a diagram:
\[
\mathcal{D}'_k(\sigma\mapsto\alpha)\colon s\stackrel{\mathcal{D}'_k(\sigma\mapsto\gamma)}\longrightarrow y \stackrel{\mathcal{D}'_k(\gamma\mapsto\alpha)}\longrightarrow x. 
\] 
In particular, if $y\in \mathcal{D}'_k(\gamma)$ and $z\in \mathcal{D}'_k(\zeta)$ are points that map by $\mathcal{D}'_k(\gamma\mapsto \alpha)$ and $\mathcal{D}'_k(\zeta\mapsto \alpha)$ onto $x$ respectively, then there is a commutative diagram 
\[
\begin{CD}
 s @>\mathcal{D}'_k(\sigma\mapsto \gamma)>> y \\
 @V \mathcal{D}'_k(\sigma\mapsto \zeta)VV @VV\mathcal{D}'_k(\gamma\mapsto \alpha)V\\
 z @>\mathcal{D}'_k(\zeta\mapsto \alpha)>> x.
\end{CD}
\] 
\end{lemma}

\begin{proof}[Proof of Lemma~\ref{l:14.2}] The point $\eta_{\alpha,k}(x)$ is in the interior of a unique non-degenerate simplex $\Delta$ of $|\mathcal{F}|$. By the definition of $|\mathcal{F}|$ the simplex $\Delta$ corresponds to an $\mathcal F$-submanifold $W$ in $\Delta\times \R^{\infty}$. Let $M$ be a compact smooth oriented manifold obtained from the fiber $W$ over $x$ by removing points in the image of the kernel of the local marking of the map $W\to \Delta$. Let $\sigma$ denote the index of the diffeomorphism type of $M$. Then the point $\eta_{\alpha,k}(x)$ belongs to the image of $\eta_{\sigma,k}$. In fact there is a unique point $s\in \mathcal{D}'_k(\sigma)$ such that $\eta_{\sigma,k}(s)=\eta_{\alpha, k}(x)$. It immediately follows that the index $\sigma$ and the point $s$ satisfy the conclusion of Lemma~\ref{l:14.2}.  
\end{proof}

\begin{proof}[Proof of Theorem~\ref{th:14.1}] For the homotopy colimits $\mathfrak{H}'_k(d)$ of $\mathcal{D}'_k$ and $\mathfrak{H}_k(d)$ of $\mathcal{D}_k$, the telescopes $\mathcal{T}'(d)$ of 
\begin{equation}\label{eq:14.1}
   \mathfrak{H}'_0(d) \longrightarrow \mathfrak{H}'_1(d) \longrightarrow \mathfrak{H}'_2(d) \longrightarrow \cdots 
\end{equation}
and $\mathcal{T}(d)$ of 
\begin{equation}\label{eq:14.2}
   \mathfrak{H}_0(d) \longrightarrow \mathfrak{H}_1(d) \longrightarrow \mathfrak{H}_2(d) \longrightarrow  \cdots
\end{equation}
are weakly homotopy equivalent to $\mathfrak{H}'(d)$ and $\mathfrak{H}(d)$ respectively. We claim that the sequences of maps (\ref{eq:14.1}) and (\ref{eq:14.2}) fit into a homotopy commutative diagram
\[
\begin{CD}
   \mathfrak{H}'_0(d) @>>> \mathfrak{H}'_1(d) @>>> \mathfrak{H}'_2(d) @>>> \cdots \\
  @VVV @VVV @VVV @. \\
   \mathfrak{H}_0(d) @>>> \mathfrak{H}_1(d) @>>> \mathfrak{H}_2(d) @>>> \cdots,
\end{CD}
\]
where each vertical map is a weak homotopy equivalence. 

\begin{lemma}\label{l:14.1} For each $k$, there is a weak homotopy equivalence
\[
    v_k: \mathfrak{H}'_k(d)\longrightarrow \mathfrak{H}_k(d).
\]
\end{lemma}

\begin{proof}[Proof of Lemma~\ref{l:14.1}] 
To begin with we arbitrarily choose weak homotopy equivalences $i_{\beta}: \mathcal{D}'_k(\beta)\to \mathcal{D}_k(\beta)$ for each $\beta\in \I_k\setminus \I_{k-1}$. Next, we assume that for some positive integer $j< k$, 

\begin{itemize}
\item[$\mathbf{A}_{j-1}\colon$] the maps $i_{\beta}$ can be constructed for all indices $\beta\in \I_{k}\setminus \I_{j}$ so that all the diagrams over arrows in $\I_{k}\setminus \I_{j}$ commute. 
\end{itemize}

We will establish the statement $\mathbf{A}_j$, and then, by induction, conclude that $\mathbf{A}_k$ holds true. 

Let $\alpha$ be an index in $\I_{j}\setminus \I_{j-1}$. We claim that there is a weak homotopy equivalence $i_{\alpha}: \mathcal{D}'_k(\alpha)\to \mathcal{D}_k(\alpha)$ such that the diagrams 
\[
\begin{CD}
   \mathcal{D}'_k(\beta) @>\mathcal{D}'_k(\beta\mapsto \alpha)>> \mathcal{D}'_k(\alpha) \\
   @V i_{\beta} VV @Vi_{\alpha} VV \\
   \mathcal{D}_k(\beta) @>\mathcal{D}_k(\beta\mapsto\alpha)>> \mathcal{D}_k(\alpha),
\end{CD}
\] 
are strictly commutative for each arrow $\beta\mapsto \alpha$ in $\I_k\setminus \I_{j-1}$. Indeed, if $x\in \mathcal{D}'_k(\alpha)$ is in the image of $\mathcal{D}'_k(\beta\mapsto\alpha)$ for some $\beta\in \I_k\setminus \I_{j}$, then by Lemma~\ref{l:14.2} there is a unique index $\sigma$ and a point $s\in \mathcal{D}'_k(\sigma)$ that satisfies the statement of Lemma~\ref{l:14.2}. We define $i_{\alpha}(x)$ to be the point $\mathcal{D}_k(\sigma\mapsto \alpha)\circ i_{\sigma}(s)$.

We recall that the space $\mathcal{D}'_k(\alpha)$ is a CW complex weakly homotopy equivalent to $\BDiff M_{\alpha}$ and there is a fiber bundle $E_{\alpha}$ over $\mathcal{D}'_k(\alpha)$ with fiber $M_{\alpha}$ associated with the principal $\Diff M_{\alpha}$ bundle. The map $i_{\alpha}$ has been constructed on a CW subcomplex $\mathop{\mathrm{Im}} D'_k(\beta\mapsto \alpha)$ of $\mathcal{D}'_k(\alpha)$ and classifies the restriction of $E_{\alpha}$. Consequently $i_{\alpha}$ can be extended to a map of $\mathcal{D}'_k(\alpha)$ so that $i_{\alpha}$ is a weak homotopy equivalence. 

Now by induction we may choose maps $i_{\gamma}$ for all $\gamma\in \I_k$ so that all diagrams over arrows in $\I_k$ are strictly commutative. Finally the maps $i_{\alpha}$ define a map $v_k$ of homotopy colimits and $v_k$ is a weak homotopy equivalence. 
 \end{proof}

\begin{lemma}\label{l:14.3}
The diagram preceding Lemma~\ref{l:14.1} with vertical maps given by the maps $v_k$ is homotopy commutative. 
\end{lemma}
\begin{proof} Let us recall that in order to construct a map $v_k$ we first arbitrarily choose maps $i_{\beta}: \mathcal{D}'_k(\beta)\to \mathcal{D}_k(\beta)$ for each $\beta\in \I_k\setminus \I_{k-1}$ and then define $v_k$ so that it extends the chosen maps $i_{\beta}$. Let us observe now that any homotopy of the chosen maps $i_{\beta}$ extends to a homotopy of $v_k$. Let us show that the diagram 
\[
\begin{CD}
   \mathfrak{H}'_k(d) @>>> \mathfrak{H}'_{k+1}(d)  \\
  @VVV @VVV  \\
   \mathfrak{H}_k(d) @>>> \mathfrak{H}_{k+1}(d),
\end{CD}
\]
is homotopy commutative. The map $v_{k+1}$ restricted to the spaces $\mathcal{D}'_k(\beta)$, with $\beta\in \I_k\setminus\I_{k-1}$, is homotopic to the disjoint union of maps $i_{\beta}$ chosen in the definition of $v_k$. By modifying $v_k$ by homotopy, we may assume that the restriction of $v_{k+1}$ to the spaces $\mathcal{D}'_k(\beta)$, with $\beta\in \I_k\setminus\I_{k-1}$ coincides with the maps $i_{\beta}$. Then we may use induction to construct a homotopy of $i_{\alpha}$ for each $\alpha\in \I_k$ to the restriction of $v_{k+1}$ to $\mathcal{D}'_k(\alpha)$. This produces a homotopy of $v_k$ to the restriction of $v_{k+1}$ to $\mathfrak{H}'_k(d)$.   
\end{proof}

The existence of maps $v_j$ that fit the above homotopy commutative diagram implies the existence of a map 
\[
   h\colon \mathcal{T}'(d)\longrightarrow \mathcal{T}(d)
\] 
of the telescopes of (\ref{eq:14.1}) and (\ref{eq:14.2}). We claim that $h$ is a weak homotopy equivalence. Indeed, each map $f$ of a sphere into $\mathcal{T}(d)$ lifts to a map into $\mathfrak{H}_j(d)$ for some finite $j$. Since $v_j$ is a weak homotopy equivalence, the latter map lifts to a map into $\mathfrak{H}'_j(d)$ and therefore determines an element in a homotopy group of $\mathcal{T}'(d)$ that maps to $[f]\in \pi_*(\mathcal{T}(d))$ by $h_*$. By a similar argument it is easily verified that $h_*$ is injective. Thus $\mathcal{T}'(d)$ is weakly homotopy equivalent to $\mathcal{T}(d)$. On the other hand, since $\mathfrak{H}'(d)$ is a CW complex, the telescope $\mathcal{T}'(d)$ is homotopy equivalent to $\mathfrak{H}'(d)$. Thus, there is a sequence of weak homotopy equivalences:
\[
   \mathfrak{H}'(d)\longrightarrow \mathcal{T}'(d)\longrightarrow \mathcal{T}(d)\longrightarrow \mathfrak{H}(d) 
\] 
and, in particular, $\mathfrak{H}'(d)$ is weakly homotopy equivalent to $\mathfrak{H}(d)$.

\end{proof}

\begin{proof}[Proof of Theorem~\ref{th:2.2}] The existence of the homotopy colimit follows from Theorem~\ref{th:5.1}. On the other hand our argument in the proof of Theorem~\ref{th:14.1} does not use any property of the chosen model of $\mathfrak{H}(d)$. Thus the argument in the proof of Theorem~\ref{th:14.1} implies that for any choices in the definition  of the homotopy colimit of classifying spaces of diffeomorphism groups of manifolds, for the resulting space $\tilde{\mathfrak{H}}(d)$, there is a weak homotopy equivalence $\mathfrak{H}'(d)\to \tilde{\mathfrak{H}}(d)$. 
\end{proof}

\section{Eliashberg-Mishachev Theorem}\label{s:5}

To determine the set of stable characteristic classes we will use the Eliashberg-Mishachev Theorem on wrinkled maps \cite{EM}, which is an h-principle type theorem.  

The {\it standard wrinkle} of type $(n,q,s)$ is defined to be a map 
\[
    w(n, q, s)\colon \R^{q-1}\times \R^{n-q}\times \R^1 \longrightarrow \R^{q-1}\times \R,
\]
parametrized by $n\ge q\ge 1$ and $s=0, ..., \left\lceil  (n-q)/2\right\rceil$, and given by 
\[
    (y, x, z) \mapsto (y, z^3+ 3(|y|^2-1)z - \sum_{i=1}^{s} x_i^2 + \sum_{j=s+1}^{n-q}x_j^2).
\]
There is an obvious homotopy of the standard wrinkle to a submersion $\R^n\to \R^q$. In particular there is a surjective homomorphism $T\R^n\to T\R^q$ covering $w(n,q,s)$. Furthermore we may choose such a homomorphism $\mathcal{R}(dw)$ so that it coincides with $dw(n,q,s)$ outside of a small neighborhood of the disc
\[
    D:=\{\, (y,x,z)\, |\, z^2+|y|^2\le 1, x=0\, \}. 
\]
The homomorphism $\mathcal{R}(dw)$ is said to be a {\it regularized differential} of $w(n,q,s)$. 

\begin{remark} We note that a homotopy of $w(n,q,s)$ to a submersion can not be chosen to be constant outside of a small neighborhood of $D$.
\end{remark}

A smooth map $f: M\to Q$ is called {\it wrinkled} if for some disjoint open subsets $U_1, ..., U_k$, the restriction $f|U_i$ is equivalent to the restriction of a standard wrinkle in a neighborhood of $D$, and for $U=\sqcup U_i$, the restriction $f|M\setminus U$ is a submersion. In view of the regularized differentials of wrinkles, every wrinkled map $f$ is covered by a {\it regularized differential} $\mathcal{R}(df)$ of $f$. 

\begin{theorem}[Eliashberg-Mishachev \cite{EM}] Given a surjective homomorphism $F: TM\to TQ$ covering a continuous map $f: M\to Q$, there is a homotopy of $f$ to a wrinkled map such that $\mathcal{R}(df)$ is homotopic to $F$. 
\end{theorem}

A wrinkled map $f$ into a manifold of dimension $>1$ is not a fold map, but it is easy to see that there is a homotopy of $f$ to such a map $\tilde f$. Furthermore, the map $\tilde f$ can be chosen so that the normal bundle of the singular set of $\tilde f$ in the source manifold is trivial. In fact, the map $\tilde f$ can be chosen so that it admits a structure of a marked fold map.

\section{The space $\Omega^{\infty}\mathbf{hV}$}\label{s:16}

In this section we recall the definition of the space $\Omega^{\infty}\mathbf{hV}$. Our definition is different from but equivalent to that given in \cite{MW}.

Let $\xi_n: E_n\to {\BSO}_n$ denote the universal orientable vector bundle of dimension $n$. Then 
\[
    \pi_n: \mathop\mathrm{Sur}(\varepsilon^{n+d}, E_n) \longrightarrow {\BSO}_n
\] 
is a vector bundle over $\BSO_n$ whose fiber over a point $b\in \BSO_n$ is a vector space of surjective homomorphisms of $\R^{n+d}$ onto the fiber $E_n|b$ of $\xi$ over $b$. Let $\mathop\mathrm{Th}\pi_n^*\xi_n$ denote the Thom space of $\pi^*\xi_n$. Then the space $\Omega^{\infty}\mathbf{hV}$ is defined to be the colimit of spaces
\[
    \Omega^{\infty}\mathbf{hV}\colon=\mathop{\mathrm{colim}}_{n\to \infty} \Omega^{n+d}\mathop\mathrm{Th}\pi_n^*\xi_n,
\] 
where the colimit is taken with respect to the obvious inclusions 
\[
    \Omega^{n+d}\mathop\mathrm{Th}\pi_n^*\xi_n \longrightarrow \Omega^{n+d}\Omega S\mathop\mathrm{Th}\pi_n^*\xi_n \longrightarrow  \Omega^{n+d+1}\mathop\mathrm{Th}\pi_{n+1}^*\xi_{n+1}.  
\]

\begin{example} If $d=0$, then the space $\mathop\mathrm{Sur}(\varepsilon^{n+d}, E_n)$ is the space of $n$ frames in $E_n$. Consequently, the space 
\[
    \mathop{\mathrm{colim}}_{n\to \infty} \mathop\mathrm{Sur}(\varepsilon^{n}, E_n)
\]
is contractible. Clearly $\Omega^{\infty}\mathbf{hV}\approx\Omega^{\infty}S^{\infty}$ in this case. 
\end{example}

The space $\Omega^{\infty}\mathbf{hV}$ plays an important role in the study of characteristic classes of manifold bundles. It is the classifying space for the cobordism group of \emph{stable formal submersions}, i.e., the cobordism group of sections 
\[
  \mathop{\mathrm{Hom}}(TM\oplus \varepsilon, TN\oplus\varepsilon)\longrightarrow M,
\]    
where $\varepsilon$ is a trivial line bundle over an arbitrary space and $TX$ denotes the tangent bundle of a smooth manifold $X$. 

\begin{remark} A \emph{formal solution} of a differential relation $\mathcal{R}\subset J^k(M,N)$ is defined to be a section of the $k$-jet bundle 
$J^k(M,N)\longrightarrow M$
with image in $\mathcal{R}$. A \emph{genuine solution} of $\mathcal{R}$ is a smooth map $f: M\to N$ whose $k$-jet extension $j^kf$ is a formal solution of $\mathcal{R}$. The \emph{homotopy principle}, or h-principle, asserts that the space of formal solutions of $\mathcal{R}$ is weakly homotopy equivalent to the space of genuine solutions of $\mathcal{R}$. 

In \cite{Sa} we introduce the notion of a \emph{stable formal solution} of a differential relation $\mathcal{R}$ with sufficiently many symmetries and formulated a (co)bordism version of the h-principle, or \emph{b-principle}. The cobordism version of the $k$-jet extension map, transforms a formal solution to $N$ into a continuous map $N\mapsto \Omega^{\infty}\mathbf{B}_{\mathcal{R}}$ into an appropriate infinite loop space. Under mild conditions\footnote{This conditions are not satisfied for submersions.} on the differential relation, we showed that cobordism groups of stable formal solutions into $N$ are isomorphic to homotopy classes of maps $[N, \Omega^{\infty}\mathbf{B}_{\mathcal{R}}]$ of $N$. 

In the case of submersions the infinite loop space $\Omega^{\infty}\mathbf{B}_{\mathcal{R}}$ coincides with $\Omega^{\infty}\mathbf{hV}$. 
\end{remark}

Let us recall the construction that associates to a submersion onto $N$
a continuous map $N\to \Omega^{\infty}\mathbf{hV}$. Let $f: M\to N$ be a submersion and $df: TM\to TN$ the differential of $f$. For sufficiently big positive integer $t$, there is an embedding $i$ of $M$ into $S^{t}$. We will denote the normal bundle of $i$ by $\nu(M)$. 
Then 
\[
     df\oplus \mathop{\mathrm{id}}\colon TM\oplus \nu(M)\longrightarrow f^*TN\oplus \nu(M)
\] 
is a surjective homomorphism of fiber bundles. We observe that $TM\oplus \nu(M)$ is canonically isomorphic to a trivial vector bundle over $M$, while $f^*TN\oplus \nu(M)$ is isomorphic to the normal bundle $\nu$ of the embedding $i\times f: M\to S^t\times N$. Consequently, a map $M\to \BSO_n$ classifying $\nu$ lifts to a map $M\to \mathop{\mathrm{Sur}}(\varepsilon^{n+d}, E_n)$. Now the Pontrjagin-Thom construction yields a map $N\to \Omega^{\infty}\mathbf{hV}$.  

\begin{remark} For each manifold $M_{\alpha}$, there is a canonical inclusion 
\[
i_{\alpha}:\BDiff M_{\alpha}\longrightarrow \mathfrak{H}(d).
\] 
On the other hand, the Pontrjagin-Thom construction yields maps 
\[
j_{\beta}: \BDiff M_{\beta}\longrightarrow \Omega^{\infty}\mathbf{hV},
\] 
one for each oriented closed manifold $M_{\beta}$. We note that in general the maps $j_{\beta}$ do not extend to a map $\mathfrak{H}(d)\to \Omega^{\infty}\mathbf{hV}$. For example, suppose that $d=2$. Then the Pontrjagin-Thom construction maps $\BDiff S^2$ and $\BDiff T^2$ to two different components in $\Omega^{\infty}\mathbf{hV}$ since the Euler characteristics of $S^2$ and $T^2$ are different (e.g., see \cite[page 766]{GMT}). On the other hand the space $\mathfrak{H}(2)$ is connected. 
\end{remark}

\section{Tautological classes}\label{sec:17}

Each class $c$ in $H^*(\Omega^{\infty}\mathbf{hV})$ determines a so called \emph{tautological} characteristic class of manifold bundles, also denoted by $c$. Namely, the class $c$ associates to a manifold bundle $\xi$ over $N$ a cohomology class $c(\xi)$ given by the pullback of $c$ with respect to the continuous map
\[
    N\longrightarrow \Omega^{\infty}\mathbf{hV}
\]
associated with $\xi$ (see section~\ref{s:16}). 

The class $c(\xi)$ can be computed using the Pontrjagin-Thom construction. 
To simplify the description let us assume that $c$ corresponds to an additive generator $P$ of a cohomology group of $\BSO_d$. Let $\xi: M\to N$ be a smooth manifold bundle. It determines a so called vertical tangent bundle $\eta$ of the manifold bundle which is a subbundle of $TM$ of planes tangent to the fibers of $f$. Then $P(\eta)$ is a cohomology class in $M$. The Pontrjagin-Thom construction applied to the map $\xi: M\to N$ yields a map
\begin{equation}\label{e:7.1}
       S^{t+d}\wedge N_+\longrightarrow \mathop{\mathrm{Th}}(\nu). 
\end{equation}
where $\nu$ is the normal bundle of $M$ in $S^{t+d}\wedge N_+$ and $\mathop{\mathrm{Th}}(\nu)$ is the Thom space of $\nu$. The composition 
\[
     H^*(M)\longrightarrow H^{*+t}(\mathop{\mathrm{Th}}(\nu)) \longrightarrow H^{*+t}(S^{t+d}\wedge N_+) \longrightarrow H^{*-d}(N)
\]
of the Thom isomorphism, the homomorphism induced by the map (\ref{e:7.1}), and again the Thom isomorphism 
takes the class $P(\eta)$ onto the desired class $c(\xi)$. \\

\section{Stable characteristic classes}\label{se:16}

We recall that the space $\sqcup \BDiff M$ is weakly homotopy equivalent to the CW complex $\D_0$ and therefore the map 
\[
     \D_0\longrightarrow \sqcup \BDiff M
\]
induces an isomorphism of cohomology groups. In particular we may identify characteristic classes of oriented manifold bundles with cohomology classes of $\D_0$. 

The path components of $\D_0$ are indexed over diffeomorphism types of oriented closed manifolds $M$ of dimension $d$. Over the component $\D_0(M)$ associated with $M$ there is a manifold bundle $\pi_M: \mathbf{W}(M)\to \D_0(M)$ with fiber $M$. We recall that the fiber $M_y$ over each point $y$ in $\D_0(M)$ is embedded into $S^{\infty}$. Let 
\[
    \varphi\colon \D_0 \longrightarrow \m.
\]
be the map given by the Pontrjagin-Thom construction applied to the disjoint union of manifold bundles $\pi_M$. In terms of section~\ref{s:16} the map $\varphi$ is associated with the disjoint union of maps 
\[
     \mathbf{W}(M)\longrightarrow \mathop{\mathrm{Sur}}(\varepsilon^{\infty+d}, E_{\infty})
\] 
that takes a point $x\in M_y$, with $y=\pi_M(x)$, to the surjective homomorphism given by the Whitney sum of the identity homomorphism $\id_{\R^d}$ and the projection of  $T_xS^{\infty}=\varepsilon^{\infty}$ to the normal plane of $M_y$ in $S^{\infty}$ at $x$. 

From the definition of the homotopy colimit $\mathfrak{H}(d)$ it follows that the path components of $\mathfrak{H}(d)$ are in bijective correspondence with cobordism classes of oriented closed manifolds of dimension $d$. Furthermore any two path components of $\mathfrak{H}(d)$ are weakly homotopy equivalent. Consequently the ring of stable characteristic classes is given by the product of rings $R_{\alpha}$ indexed by cobordism classes $\alpha$ of oriented closed manifolds of dimension $d$; for each $\alpha$, the term $R_{\alpha}$ is the ring of stable characteristic classes of oriented manifold bundles whose fibers are in the cobordism class $\alpha$. 

Let $\Omega^{\infty}_0\mathbf{hV}$ denote the path component of loops in $\Omega^{\infty}\mathbf{hV}$ homotopic to the trivial loop. All path components of $\Omega^{\infty}\mathbf{hV}$ are homotopy equivalent, and there is a projection
\[
\psi\colon \m \longrightarrow \Omega^{\infty}_0\mathbf{hV}
\]
whose restriction to each component is a homotopy equivalence. It is easily verified that the pullback $\varphi^*(\psi^*c)$ of any class $c$ in $H^*(\Omega^{\infty}_0\mathbf{hV})$ is a characteristic class.

\begin{theorem}\label{th:6.3} Suppose that $d$ is even. Then each stable characteristic class in $R_{\alpha}$ is of the form $\varphi^*(\psi^*c)$ for some class $c$ in $H^*(\Omega^{\infty}_0\mathbf{hV})$. 
\end{theorem}
\begin{proof} 
A manifold bundle $f: W\to N$ over a closed oriented path connected manifold $N$ of dimension $n$ with fiber an oriented closed manifold $M$ of dimension $d$ is classified by a map $N\to \BDiff M$, whose composition with the inclusion $\BDiff M \hookrightarrow \mathfrak{H}(d)$ is a map 
\[
\tau(f)\colon N\longrightarrow \mathfrak{H}(d). 
\]
In particular, each manifold bundle $f$ determines a cobordism class $[\tau(f)]$ in $\Omega_n(\mathfrak{H}(d))$. There is an isomorphism of groups
\begin{equation}\label{eq:18.1}
\Omega_n(\mathfrak{H}(d))\otimes\Bbbk \cong \sum_{p+q=n} H_p(\mathfrak{H}(d); \Omega_q\otimes \Bbbk),
\end{equation}
where $\Omega_q$ is the $q$-th cobordism group of a point \cite{CF}. We say that a manifold bundle $f$ is \emph{prime} if the element $[\tau(f)]$ is taken by the isomorphism (\ref{eq:18.1}) onto an element with only one non-trivial component in $H_n(\mathfrak{H}(d); \Bbbk)$. We note that each element in $H_n(\mathfrak{H}(d); \Bbbk)$ is the image under the isomorphism (\ref{eq:18.1}) of an element $[\tau(f)]$ for some prime manifold bundle $f$. 

\begin{lemma}\label{l:10} Let $c$ be a stable characteristic class of degree $n$ in $R_{\alpha}$. If $f_i: W_i\to N_i$, with $i=0,1$,
are two prime manifold bundles, whose fibers are in the cobordism class of $\alpha$, over closed oriented path connected manifolds $N_i$ of dimension $n$ such that $\varphi_*([f_0])$ and $\varphi_*([f_1])$ are the same in the group $H_n(\m)$, then $c(f_0)=c(f_1)$. 
\end{lemma}
\begin{proof}
There is a commutative diagram of homomorphisms of groups
\[
\begin{CD}
 \Omega_n(\mathbf{T}_0)\otimes\Bbbk @>\varphi_*>> \Omega_n(\Omega^{\infty}\mathbf{hV})\otimes\Bbbk \\
 @VVV @VVV \\
 \sum_{p+q=n} H_p(\mathbf{T}_0; \Omega_q\otimes \Bbbk) @>>> \sum_{p+q=n} H_p(\Omega^{\infty}\mathbf{hV}; \Omega_q\otimes \Bbbk) \\
 @VVV @VVV\\
 H_n(\mathbf{T}_0) @>\varphi_*>> H_n(\Omega^{\infty}\mathbf{hV}),  
\end{CD}
\]
where the lower vertical homomorphisms are projections onto the $n$-th factors. From the top half of the diagram it follows that each of the cobordism classes $\varphi_*([f_0])$ and $\varphi_*([f_1])$ maps to an element in 
$\sum_{p+q=n} H_p(\Omega^{\infty}\mathbf{hV}; \Omega_q\otimes \Bbbk)$ with only one non-trivial component and that component is in $H_n(\Omega^{\infty}\mathbf{hV})$. Consequently,  $\varphi_*([f_0])=\varphi_*([f_1])$ in $\Omega_n(\Omega^{\infty}\mathbf{hV})\otimes\Bbbk$. 

Thus, without loss of generality, we may assume that there is a map $f: W\to N\times [0,1]$ with $\partial W=W_0\sqcup -W_1$ and $\partial N= N_0\sqcup -N_1$ which restricts to 
\[
f_0: W_0\to N_0\times \{0\} \qquad  \mathrm{and} \qquad f_1: W_1\to N_1\times \{1\}
\] 
and such that there is an epimorphism $TW\oplus \varepsilon\to TN\oplus \varepsilon$, where $\varepsilon$ is a trivial line bundle, covering $f$ that restricts to $df_0\oplus \id_{\varepsilon}$ and $df_1\oplus \id_{\varepsilon}$ over the boundary. Then, by the Eliashberg-Mishachev Theorem, the map $f\times \id_{S^1}$ is homotopic to a wrinkled map
\[
     g: W\times S^1\longrightarrow N\times S^1
\]
 under homotopy which is trivial near the boundary. Furthermore the map $g$ is homotopic to a marked fold map $\tilde g$ under homotopy trivial near the boundary. 
 
Let $pt$ be a point in $S^1$. We may modify $g$ by a homotopy so that the pull back of $g$ with respect to the inclusion 
\[
  N\stackrel{=}\longrightarrow N\times \{pt\}\stackrel{\subset}\longrightarrow N\times S^1
\] 
is a marked fold map which restricts to families $f_0$ and $f_1$ over the boundary $\partial N$.  Hence $c(f_0)=c(f_1)$.
\end{proof}

The finite dimensional vector space $H_n(\m)$ decomposes into the sum of two vector subspaces
\[
    H_n(\m) \stackrel{\cong}\longrightarrow \mathop{\mathrm{Im}}(\varphi_*) \oplus \mathop{\mathrm{Coker}}(\varphi_*).
\]
By Lemma~\ref{l:10} there is a class $\tilde c \in \mathop{\mathrm{Hom}}(H_*(\m), \Bbbk)$ defined trivially on $\mathop{\mathrm{Coker}}(\varphi_*)$ such that it is compatible with $\tilde c$ on $\mathop{\mathrm{Im}}(\varphi_*)$. 
\end{proof}

\addcontentsline{toc}{part}{Stability condition for characteristic classes}

%\section{Stability condition for coverings}
%
%In this section we consider manifold bundles with fiber given by an oriented compact %manifold of dimension $0$. In other words, we consider coverings $f: X\to Y$ each fiber %$f^{-1}(y)$ of which is a finite set of oriented points and the orientation of the points %in the fiber $f^{-1}(y)$ changes continuously as $y$ ranges over points in $Y$. For %example, the projection $\mathbf{1}_Y: S^0\times Y\to Y$ is an oriented covering each %fiber of which consists of one positive and one negative point. If $\xi_1: X_1\to Y$ and %$\xi_2: X_2\to Y$ are two oriented coverings, then 
%\[
%   \xi_1\sqcup \xi_2: X_1\sqcup X_2 \longrightarrow Y
%\] 
%is also an oriented covering. 
% 
%Let $c$ be a characteristic class of oriented coverings. Then $c$ is \emph{stable} if 
%\[
%    c(\xi)=c(\xi\sqcup \mathbf{1}_Y)
%\]  
%for each oriented covering $\xi$ over each topological space $Y$. 

\section{Geometric interpretation of stability conditions of low order}\label{s:2}

In this section we give a geometric interpretation of stability and strong stability conditions of low order. To begin with, let us consider the case $d=2$, where the stability of characteristic classes has been studied in relation to moduli spaces of Riemann surfaces and the standard Mumford conjecture in particular. 
%by Madsen and Weiss 
%turned out to be polynomials in terms of even Miller-Morita-Mumford classes. 

\subsection{Stability condition of order $1$}

Let $\xi_g$ be a surface bundle over a space $X$. %with fiber given by a surface $F_g$ of genus $g$. 
Suppose that $\xi_g$ admits a pair of disjoint sections $s_0$ and $s_1$. Then each fiber $F$ of $\xi_g$ is penetrated by $s_0$ and $s_1$ at two points, say $p$ and $q$. If tubular neighborhoods of $s_0(X)$ and $s_1(X)$ are trivial disc bundles over $X$ with a fixed trivialization on each, then we can attach a $1$-handle to each fiber $F$ along $p$ and $q$ in a consistent obvious way to obtain a new manifold bundle $\xi_{g+1}$ over $X$, which is a {\it $1$-handle stabilization of $\xi_g$} in the sense that if the fiber of $\xi_g$ is a surface $F_g$ of genus $g$, then the fiber of $\xi_{g+1}$ is a surface $F_{g+1}$ of genus $g+1$. A {\it $0$-handle stabilization} of $\xi_g$ is defined to be the fiberwise disjoint union $\xi_{g+0}$ of $\xi_g$ and the trivial $S^q$-bundle over $X$. 
%The characteristic class $c$ associates to $\xi_{g+1}$ a new cohomology class 
%$c(\xi_{g+1})$ in $H^t(Q)$ which may differ from $c(\xi_g)$. 
The mentioned {\it stability condition of order $1$} on a characteristic class $c$ of surface bundles means that 
\[
 c(\xi_g)=c(\xi_{g+1})=c(\xi_{g+0}),
\]
i.e., $c$ does not change under stabilizations of surface bundles. 

If dimension $d$ is higher than $2$, then we need to consider stabilizations of fiber bundles by attaching to each fiber a handle of higher index $i$ as well. Namely, let $\xi_{\alpha}$ be a manifold bundle over a space $X$ with fiber given by an oriented manifold $M_{\alpha}$ of dimension $d$. Suppose that in each fiber $F$ there is an embedded sphere $S$ of dimension $i-1$ with embeddings chosen in consistent way so that the spheres form a manifold bundle $E\to X$ whose structure group reduces to the orthogonal group $O_i$. Suppose that each embedded sphere $S$ has a trivial normal bundle in the fiber $F$ with a fixed trivialization smoothly depending on the fiber $F$. Then we can attach a handle of index $i$ to each fiber $F$ along $S$ in a consistent obvious way to form a manifold bundle $\xi_{\beta}$ with fiber given by a manifold $M_{\beta}$.  We say that the fiber bundle $\xi_{\beta}$ is a {\it stabilization} of $\xi_{\alpha}$.

\begin{definition} A characteristic class $c$ of manifold bundles is {\it stable of order $1$} if for any manifold bundle $\xi_{\alpha}$ of dimension $d$ and any stabilization $\xi_{\beta}$ of $\xi_{\alpha}$, there is an equality of cohomology classes $c(\xi_{\alpha})=c(\xi_{\beta})$.  
\end{definition}

%\begin{remark} An inductive application of the Mayer-Vietoris theorem shows %that Definitions~\ref{d:1.1} and \ref{d:2.3} agree.  
%\end{remark}

%\begin{remark} Each {\it higher torsion invariant} \cite{Ig} is a stable characteristic %class of order $1$ on an appropriate class of manifold bundles. In particular, for each %$k>0$, the higher Miller-Morita-Mumford class $M_{2k}$ and 
%\[
%       M_{2k}(E):=tr^{E}_Q((2k)!\,ch_{4k}(T^vE))\in H^{4k}(Q; \Z)
%\]
%the higher Franz-Reidemeister torsion invariant $\tau_{2k}$ are stable characteristic %classes on the class of fiber bundles over which these invariants are well-defined, e.g., %on the class of oriented fiber bundles $E\to Q$ with fibers given by closed manifolds $F$ %such that the fundamental group $\pi_1(Q)$ acts trivially on $H_*(F; \Q)$ (see %Proposition~\ref{l:3.1}).  \end{remark}

The classes $RW_{2k}$ in Example~\ref{ex:1.1} are non-trivial, but stably trivial in the sense that $RW_{2k}(\xi)=0$ for almost all surface bundles $\xi$. As we will see this implies that the classes $RW_{2k}$ violate a higher stability condition.  

%A characteristic class $c$ is said to be {\it additive} if for any two %manifold bundles $\xi, \eta$ over a space $X$, there is an equality 
%\[ 
% c(\xi\sqcup \eta)=c(\xi)+c(\eta),
% \] 
% where $\xi\sqcup \eta$ stands for the manifold bundle over $X$ given by the %fiberwise disjoint union of $\xi$ and $\eta$.

\subsection{Stability condition of order $2$}\label{s:2.5}

We recall that a stable characteristic class of order $0$ can be associated with a cohomology class in $H^{*}(\mathfrak{H}_0)$. Such a class is of order $1$ if it extends to a class in $H^*(\mathfrak{H}_1)$. The space $\mathfrak{H}_1=\mathfrak{H}_1(d)$ has many new cycles that are not present in $\mathfrak{H}_0(d)$. The simplest ones can be described as follows. Let $S^i$ and $S^j$ be two disjoint spheres embedded into a closed oriented manifold $M_0$ with trivialized normal bundles. Let $U_i$ and $U_j$ be disjoint neighborhoods of $S^i$ and $S^j$ respectively. Then attaching handles $H_i$ and $H_j$ along $S^i$ and $S^j$ respectively results in three new manifolds; namely $M_i$ and $M_j$ obtained by attaching handles $H_i$ and $H_j$ respectively and a manifold $M_{ij}$ obtained by attaching both $H_i$ and $H_j$. Now the diagram
%\[
%\xymatrix
%{
%& \BDiff M_0                                &\BDiff M_i                       %                  & \BDiff M_{ij}                        &\BDiff M_j   \\     
%& \BDiff M_0\setminus U_i \ar[u] \ar[ur]_{} &\BDiff M_i\setminus U_j \ar[u] %\ar[ur] & \BDiff M_j\setminus U_i\ar[u]\ar[ur] &\BDiff M_0\setminus %U_j\ar[u]\ar[ulll]
%}
%\]
%
%
\[
\begin{CD}
\BDiff M_0 @<<< \BDiff M_0\setminus U_i @>>> \BDiff M_i\\
@AAA @. @AAA \\
\BDiff M_0\setminus U_j @. @. \BDiff M_i\setminus U_j \\
@VVV @. @VVV \\
\BDiff M_j @<<< \BDiff M_j\setminus U_i @>>> \BDiff M_{ij}.
\end{CD}
\]
reveals many simple cycles in $\mathfrak{H}_1$, which we call {\it special of order $2$}. Such a cycle is the image of the fundamental class of $N\times S^1$, where $N$ is a closed oriented manifold, under the map that takes points of each copy $\{pt\}\times S^1$ consecutively to 
\[
   \BDiff M_0,\quad \BDiff M_i,\quad  \BDiff M_{ij},\quad  \BDiff M_{j}
\]
and finally to $\BDiff M_0$.

Since special cycles of order $2$ are not present in $\mathfrak{H}_0$ it is reasonable to require that each characteristic class with compatible components $c_M\in H^*(\BDiff M)$ has an extension $c$ in $H^*(\mathfrak{H}_1)$ such that $c$ evaluates trivially on special cycles of order $2$. 

\begin{proposition}\label{d:1.5} A stable characteristic class of order $1$ is stable of order $2$ if it has an extension in $H^*(\mathfrak{H}_1)$ trivial on all special cycles of order $2$. 
\end{proposition}

To construct a space for stable classes of order $2$, we attach new strata to $\mathfrak{H}_1$ along the special cycles of order $2$. This however creates new special cycles in $\mathfrak{H}_2$, the ones of order $3$. The order of stability can be determined by induction. This is the purpose of our next section.

%\begin{proposition}\label{d:1.6} A stable characteristic class of order $k-1$ %is of order $k$ if it has an extension to $\mathfrak{H}_{k-1}$ trivial on %special cycles of order $k$.  
%\end{proposition}  

\section{Criteria for stability conditions}

A characteristic class $c$ associates to each manifold bundle $\xi$ a cohomology class $c(\xi)$. Our purpose is to find a condition on $c$ which is necessary and sufficient for $c$ to be stable. We will establish Theorem~\ref{th:7.1}.

\begin{definition}\label{def:20.1}
We say that two marked fold maps $\xi_0$ and $\xi_2$ to a manifold $Q$ are {\it concordant} if there is a marked fold map $\xi$ to $Q\times [0,1]$ such that $\xi|Q\times \{i\}=\xi_i$ for $i=0,1$.  
\end{definition}

Let $\psi_M: M\to X$ and $\psi_N: N\to X$ be two fiber bundles over a CW complex $X$ such that for each point $x\in X$, the fibers $M_x$ and $N_x$ of respectively $\psi_M$ and $\psi_N$ over $x$ are smooth manifolds. A map $f:M\to N$ is said to be a {\it map over $X$}, if the diagram
\[
\xymatrix
{
& M \ar[rr]^{f}\ar[dr]^{\psi_M} &       & N  \ar[dl]_{\psi_N}   \\     
&            & X& 
}
\]
commutes. Suppose that for each cell $e: D\to X$ with center $x\in D$ there are fiberwise diffeomorphisms 
\[
    M_x\times D\longrightarrow M 
\]
and
\[
   N_x\times D\longrightarrow N
\]
over $D$ such that 
\[
    M_x\times \{x\}\longrightarrow M_x \qquad \mathrm{and}\qquad N_x\times \{x\}\longrightarrow N_x
\]
are identity maps and the diagram
\[
\begin{CD}
   M_x\times D @>>> M\\
   @Vf_x\times \id_D VV @VVV\\
   N_x\times D@>>> N
\end{CD}
\]
over $D$ commutes. Then the map $f$ is said to be a {\it marked fold map over $X$} if each $f_x:=f|f^{-1}(x)$ is a marked fold map. 
 
\begin{definition} We say that two marked fold maps $\xi_i: M_i\to N$ over a space $X$, with $i=0,1$, are {\it concordant}, if there is a marked fold map $\xi: M\to N\times [0,1]$ over $X$ of an oriented compact manifold $M$ with boundary $\partial M\equiv M_0\sqcup (-M_1)$ such that $\xi^{-1}(N\times \{i\})=M_i$ and $\xi|M_i: M_i\to N\times \{i\}$ coincides with $\xi_i$ for $i=0,1$.
\end{definition}

To formulate a criteria for a stable characteristic class of order $k-1$ to be of order $k$ we need the notion of a {\it special family}. 

\begin{definition}[Special family]\label{d:8} Let $k$ be a positive integer and $X$ a CW complex. Let $g: W\to D^k$ be an elementary $k$-bordism and $\psi$ a kernel of $(W, g)$. By Lemma~\ref{l:4.1} it is associated with a manifold $M$ and an equivalence class $[w]$ of parametrizations of $\partial M$. We recall that $D^k\times M$ is canonically embedded into $W$ so that the restriction of $g$ to $D^k\times M$ is a trivial fiber bundle with projection onto the first component. Let $\id_{D^k}$ denote the trivial diffeomorphism of $D^k$. Then there is an inclusion
\[
   G\colon= \Diff (M,w) \longrightarrow \Diff (D^k\times M),
\]
\[
   g\mapsto \id_{D^k}\times g,
\]
which defines an action of $G$ on $D^k\times M$ and, furthermore, on $W$ by a trivial extension.  We note that the obtained action of $G$ on $W$ commutes with the map $g$. Let $EG\to X$ be an arbitrary principle $G$-bundle over $X$. Then $EG\times_G W$ is the total space of a family $\eta$ 
\[
\xymatrix
{
& EG\times_G W \ar[rr] \ar[dr]&       &D^k\times X \ar[dl]     \\     
&            & X & 
}
\]
over $X$ of elementary $k$-bordisms. A {\it special family} of order $k$ is defined to be a family of the form 
\[
\xymatrix
{
& X_k \ar[rr]^{\eta|X_k} \ar[dr]&       &\partial D^k\times X \ar[dl]     \\     
&            & X & 
}
\]
where $X_k$ stands for $\eta^{-1}(\partial D^k\times X)$.  
\end{definition}

\begin{example} Let $f: W\to X$ be a marked fold map into a closed oriented manifold $X$. Let $\Sigma_k$ denote the set of those parameters $x\in X$ for which the fiber $f^{-1}(x)$ has precisely $k$ singular points. Suppose that non of the fibers of $f$ has more than $k$ singular points. Then $\Sigma_k$ is a closed oriented submanifold of $X$. Furthermore the restriction of the family $f$ to the boundary of a tubular neighborhood of $\Sigma_k$ in $X$ is a special family of order $k$.  
\end{example}

\begin{theorem}\label{th:7.1} Let $c$ be a characteristic class of oriented manifold bundles of dimension $d$. If $c$ is stable of order $k+1$ and $c$ is trivial on special cycles of order $k$, then $c$ is stable of order $k+1$.  
\end{theorem}

\begin{example} In the case where $k=1$ and $N$ is an oriented closed manifold, a special family parametrized by $N\times S^0$ is given by two manifold bundles over $N$ and $-N$ such that one is a stabilization of the other. Here $-N$ stands for the manifold $N$ with orientation different from the fixed one. Similarly it is easily verified that in the case $d=2$, Theorem~\ref{th:7.1} reduces to Proposition~\ref{d:1.5}. 
\end{example}

\begin{proof}[Proof of Theorem~\ref{th:7.1}]
Let $c$ be a characteristic class which is stable of order $k$. Then $c$ can be identified with a cohomology class in $H^n(\D_k)$ for some $n$.

\begin{remark} This identification is not unique. The difference between two choices is a phantom class. However, for our argument the choice of a class in $H^n(\D_k)$ is not essential. 
\end{remark}

Suppose that $c$ is trivial on special cycles of order $k+1$. Let us show that it is stable of order $k+1$, i.e., the class $c$ lifts to a cohomology class in $H^n(\D_{k+1})$ with respect to the exact sequence  
\[
    H^{n}(\D_{k+1})\longrightarrow H^n(\D_{k})\stackrel{\partial k}\longrightarrow H^n(\D_{k+1}, \D_k).
\]
We may regard $c$ to be a an additive $\Bbbk$-valued function on $H_{n}(\D_k)$. Thus in order to show that $c$ is of order $k+1$, we need to show that it extends to an additive $\Bbbk$-valued function on $H_{n}(\D_{k+1})$ with respect to the homomorphism of vector spaces
\[
i_*\colon   H_n(\D_k)\longrightarrow H_n(\D_{k+1})
\]
induced by the inclusion $i: \D_k\subset \D_{k+1}$. To this end, it suffices to show that $c$ it trivial on the kernel of $i_*$. 

Suppose that $c$ is not trivial on the kernel of $i_*$, i.e., there is a homology class $x\in H_{n}(\D_k)$ such that $i_*(x)=0$, but $c(x)\ne 0$. 

Let us recall an isomorphism
\begin{equation}\label{eq:9.7}
\Omega_n(\D_k)\otimes\Bbbk\stackrel{\cong}\longrightarrow \sum_{p+q=n} H_p(\D_k; \Omega_q\otimes \Bbbk)
\end{equation}
where $\Omega_n(\D_i)$ stands for the oriented $n$-th bordism group of $\D_i$, while $\Omega_q$ stands for the $q$-th oriented bordism group of a point. Let  $[\alpha]\in \Omega_n(\D_i)\otimes \Bbbk$ be the element taken by the isomorphism (\ref{eq:9.7}) to the class $\tilde x$ whose first $n$ components are zero and whose last component is $x\in H_n(\D_k; \Omega_0\otimes \Bbbk)$. We observe that the  homomorphism 
\[
   \sum_{p+q=n} H_p(\D_k; \Omega_q\otimes \Bbbk)\longrightarrow \sum_{p+q=n} H_p(\D_{k+1}; \Omega_q\otimes \Bbbk)
\]
induced by the inclusion $\D_k\to \D_{k+1}$ takes $\tilde x$ to zero. Hence $i_*$ takes $[\alpha]$ to zero in $\Omega_n(\D_{k+1})\otimes \Bbbk$.

Let $\alpha$ be a map to $\D_k\subset \D_{k+1}$ representing $[\alpha]\in  \Omega_n(\D_k)\otimes \Bbbk$. Then there is a bordism $A: M\to \D_{k+1}$ of $\alpha$ to zero. Furthermore, we may assume that  $A$ is transversal to $\D_{k+1}\setminus \D_{k}$. Then the set $\Sigma\colon=A^{-1}(\D_{k+1}\setminus \D_{k})$ is an embedded oriented closed submanifold in $M$ with trivial normal bundle. Furthermore, the restriction of $A$ to the boundary $\partial U$ of a closed tubular neighborhood $U$ of $\Sigma$ is a map into $\D_{k}$ corresponding to a special family of order $k+1$. Hence, by the hypothesis of Theorem~\ref{th:7.1}, 
\[
\left\langle\ c,\ [A|\partial U]\ \right\rangle=0,
\] 
where $[A|\partial U]$ is the image in $H_n(\D_{k})$ of $(A|\partial U)_*$ of the orientation class $[\partial U]$.  
Consequently, the value 
\begin{equation}\label{eq:9.8}
\left\langle\ c,\ [A|\partial(M\setminus \mathop{\mathrm{Int}}(U))]\ \right\rangle  
\end{equation}
is the same as $\left\langle c, [\alpha] \right\rangle\ne 0$, where $\mathop{\mathrm{Int}}(U)$ stands for the interior of $U$. On the other hand, the map $A|M\setminus \mathop{\mathrm{Int}}(U)$ is a bordism of $A|\partial(M\setminus \mathop{\mathrm{Int}}(U))$ to zero in $\D_{k}$, and therefore   
\[
   [A|\partial(M\setminus \mathop{\mathrm{Int}}(U))]=0
\]
in $\Omega_n(\D_k)$, which implies that the value of the expression (\ref{eq:9.8}) is zero. Thus the class $[\alpha]$, and therefore $x$, does not exist. 

\end{proof}

\section{Madsen-Weiss homotopy colimit decompositions}\label{MW}

Homotopy colimits closely related to $\mathfrak{H}(d)$ are studied in the paper \cite{MW} by Madsen and Weiss. The authors introduce sheaves $\mathcal{W}$ and $\mathcal{W}_{\textrm{loc}}$ of genuine objects and sheaves $h\mathcal{W}$  and $h\mathcal{W}_{\textrm{loc}}$ of their respective homotopy theoretic counterparts. Homotopy theoretic counterparts are simpler and can be studied using homotopy theory. The study of sheaves of genuine objects is more technically demanding and rests on homotopy colimit decompositions 
\[
   |\mathcal{W}|= \mathop{\mathrm{hocolim}}_{T} |\mathcal{W}_T|
\]
and 
\[
   |\mathcal{W}_{\textrm{loc}}|= \mathop{\mathrm{hocolim}}_{T} |W_{\textrm{loc}, T}|.
\]
These are obtained in \cite[Chapter 5]{MW} using sheaf theoretic technique which is developed in the paper. In the case $d=2$, Madsen and Weiss substitute the classifying spaces $|h\mathcal{W}|$ and $|h\mathcal{W}_{\textrm{loc}}|$ for the two homotopy colimits in the fibration  
\begin{equation}\label{e:3.1}
\Z\times \BDiff F^{+}_{\infty}\longrightarrow \mathop{\mathrm{hocolim}}_{T} |\mathcal{W}_T| \longrightarrow \mathop{\mathrm{hocolim}}_T |\mathcal W_{\textrm{loc}, T}|
\end{equation}
and thus determine the homotopy type of $\Z\times \BDiff F^{+}_{\infty}$, which is the content of the generalized Mumford Conjecture. 

Our approach can be regarded to be a geometric version of Madsen-Weiss approach to sheaves of genuine objects. In fact, in the notation of Madsen and Weiss, our space $\mathfrak{H}(d)$ is defined to be the colimit
\begin{equation}\label{e:3.2}
    \mathfrak{H}(d)\colon=\mathop{\mathrm{hocolim}}_T (|\mathcal{W}_T| \to |\mathcal W_{\textrm{loc}, T}|)
\end{equation}
of fibers of $|\mathcal{W}_T|\to |\mathcal{W}_{\textrm{loc}, T}|$. Of course the operations of taking fibers and homotopy colimits do not commute. Because of that the fiber of homotopy colimits in (\ref{e:3.1}) is not quite the same as the homotopy colimit of fibers in (\ref{e:3.2}).  

\begin{remark}
The existence part of Theorem~\ref{th:2.2} can be deduced from the existence of the homotopy colimit of fibers in (\ref{e:3.2}); our construction however is somewhat more explicit. The uniqueness part of Theorem~\ref{th:2.2} does not readily follow from the results in \cite{MW}. 
\end{remark}

\section{Final remark}

Similarly stable stable characteristic classes can be defined for non-oriented manifold bundles. The assumption on the coefficients of cohomology groups can also be lifted. 

\begin{example} The {\it orientation class} of non-oriented manifold bundles is a stable characteristic class with coefficients in $\Z_2$. To a non-oriented manifold bundle $\xi$ over a base space $X$, the orientation class associates a cohomology class $c(\xi)$ in $H^1(X; \Z_2)$ defined so that for a homology class in $H_1(X, \Z_2)$ represented by a closed connected loop $l: S^1\to X$, the value of $c(\xi)$ on $[l]$ is non-trivial if and only if the manifold bundle pulled back to $S^1$ by means of $l$ is non-orientable.     
\end{example}

\end{document}